\documentclass[a4paper,11pt]{amsart}
\usepackage{amsfonts,amsmath,amsthm,amssymb,color}
\usepackage[all]{xy}
\usepackage[pagebackref, colorlinks=true, linkcolor=blue, citecolor=blue]{hyperref}
\usepackage{faktor}

\DeclareMathOperator*{\colim}{colim}

\numberwithin{equation}{section}

\begin{document}

\newtheorem{theorem}{Theorem}[section]
\newtheorem{thm}[theorem]{Theorem}
\newtheorem{lemma}[theorem]{Lemma}
\newtheorem{proposition}[theorem]{Proposition}
\newtheorem{corollary}[theorem]{Corollary}

\theoremstyle{definition}
\newtheorem{definition}[theorem]{Definition}
\newtheorem{example}[theorem]{Example}

\theoremstyle{remark}
\newtheorem{remark}[theorem]{Remark}

\newenvironment{magarray}[1]
{\renewcommand\arraystretch{#1}}
{\renewcommand\arraystretch{1}}

\newcommand{\quot}[2]{
{\lower-.2ex \hbox{$#1$}}{\kern -0.2ex /}
{\kern -0.5ex \lower.6ex\hbox{$#2$}}}

\newcommand{\mapor}[1]{\smash{\mathop{\longrightarrow}\limits^{#1}}}
\newcommand{\mapin}[1]{\smash{\mathop{\hookrightarrow}\limits^{#1}}}
\newcommand{\mapver}[1]{\Big\downarrow
\rlap{$\vcenter{\hbox{$\scriptstyle#1$}}$}}
\newcommand{\liminv}{\smash{\mathop{\lim}\limits_{\leftarrow}\,}}

\newcommand{\specif}[2]{\left\{#1\,\left|\, #2\right. \,\right\}}

\renewcommand{\bar}{\overline}
\newcommand{\de}{\partial}
\newcommand{\debar}{{\overline{\partial}}}
\newcommand{\per}{\!\cdot\!}
\newcommand{\Oh}{\mathcal{O}}
\newcommand{\sA}{\mathcal{A}}
\newcommand{\sB}{\mathcal{B}}
\newcommand{\sC}{\mathcal{C}}
\newcommand{\sE}{\mathcal{E}}
\newcommand{\sF}{\mathcal{F}}
\newcommand{\sG}{\mathcal{G}}
\newcommand{\sW}{\mathcal{W}}
\newcommand{\sH}{\mathcal{H}}
\newcommand{\sI}{\mathcal{I}}
\newcommand{\sJ}{\mathcal{J}}
\newcommand{\sL}{\mathcal{L}}
\newcommand{\sM}{\mathcal{M}}
\newcommand{\sP}{\mathcal{P}}
\newcommand{\sU}{\mathcal{U}}
\newcommand{\sQ}{\mathcal{Q}}
\newcommand{\sR}{\mathcal{R}}
\newcommand{\sV}{\mathcal{V}}
\newcommand{\sX}{\mathcal{X}}
\newcommand{\sY}{\mathcal{Y}}
\newcommand{\sN}{\mathcal{N}}
\newcommand{\sS}{\mathcal{S}}

\newcommand{\Ni}{\mathcal{N}_I}
\newcommand{\Nj}{\mathcal{N}_J}
\newcommand{\bM}{\mathbf{M}}
\newcommand{\bC}{\mathbf{C}}
\newcommand{\LI}{\mathbf{L}_{\sN}}
\newcommand{\pal}{\boldsymbol{\cdot}}

\newcommand{\Aut}{\operatorname{Aut}}
\newcommand{\Mor}{\operatorname{Mor}}
\newcommand{\Def}{\operatorname{Def}}
\newcommand{\Hom}{\operatorname{Hom}}
\newcommand{\Hilb}{\operatorname{Hilb}}
\newcommand{\HOM}{\operatorname{\mathcal H}\!\!om}
\newcommand{\DER}{\operatorname{\mathcal D}\!er}
\newcommand{\Spec}{\operatorname{Spec}}
\newcommand{\Sh}{\operatorname{Sh}}
\newcommand{\Der}{\operatorname{Der}}
\newcommand{\Tor}{{\operatorname{Tor}}}
\newcommand{\Ext}{{\operatorname{Ext}}}
\newcommand{\End}{{\operatorname{End}}}
\newcommand{\END}{\operatorname{\mathcal E}\!\!nd}
\newcommand{\Image}{\operatorname{Im}}
\newcommand{\Id}{\operatorname{Id}}
\newcommand{\coker}{\operatorname{coker}}
\newcommand{\res}{\operatorname{res}}
\newcommand{\tot}{\operatorname{tot}}
\newcommand{\cone}{\operatorname{cone}}
\newcommand{\cocone}{\operatorname{cocone}}
\newcommand{\mA}{\mathfrak{m}_{A}}
\newcommand{\g}{\mathfrak{g}}

\newcommand{\somdir}[2]{\hbox{$\mathrel
{\smash{\mathop{\mathop \bigoplus\limits_{#1}}
\limits^{#2}}}$}}
\newcommand{\tensor}[2]{\hbox{$\mathrel
{\smash{\mathop{\mathop \bigotimes\limits_{#1}}
^{#2}}}$}}
\newcommand{\symm}[2]{\hbox{$\mathrel
{\smash{\mathop{\mathop \bigodot\limits_{#1}}
^{#2}}}$}}
\newcommand{\external}[2]{\hbox{$\mathrel
{\smash{\mathop{\mathop \bigwedge\limits_{#1}}
^{\!#2}}}$}}

\renewcommand{\Hat}[1]{\widehat{#1}}
\newcommand{\dual}{^{\vee}}
\newcommand{\desude}[2]{\dfrac{\de #1}{\de #2}}

\newcommand{\A}{\mathbb{A}}
\newcommand{\N}{\mathbb{N}}
\newcommand{\R}{\mathbb{R}}
\newcommand{\Z}{\mathbb{Z}}
\renewcommand{\H}{\mathbb{H}}
\renewcommand{\L}{\mathbb{L}}
\newcommand{\proj}{\mathbb{P}}
\newcommand{\K}{\mathbb{K}\,}
\newcommand\C{\mathbb{C}}
\newcommand\Del{\operatorname{Del}}
\newcommand\D{\operatorname{D}}
\newcommand\Tot{\operatorname{Tot}}
\newcommand\Grpd{\mbox{\bf Grpd}}

\newcommand\é{\'e}
\newcommand\è{\`e}
\newcommand\à{\`a}
\newcommand\ì{\`i}
\newcommand\ù{\`u}
\newcommand\ò{\`o }


\newcommand{\rh}{\rightarrow}
\newcommand{\contr}{{\mspace{1mu}\lrcorner\mspace{1.5mu}}}

\newcommand{\bi}{\boldsymbol{i}}
\newcommand{\bl}{\boldsymbol{l}}

\newcommand{\MC}{\operatorname{MC}}
\newcommand{\TW}{\operatorname{TW}}
\newcommand{\DGMod}{\operatorname{DGMod}}
\newcommand{\Ch}{\operatorname{Ch}}
\newcommand{\Ho}{\operatorname{Ho}}
\newcommand{\PM}{\mathbf{Mod}}
\newcommand{\Mod}{\mathbf{Mod}}
\newcommand{\DGPM}{\mathbf{Mod}^{\ast}}
\newcommand{\QCoh}{\mathbf{QCoh}}
\newcommand{\DGSch}{\mathbf{DGSch}}
\newcommand{\DGAff}{\mathbf{DGAff}}

\newcommand{\Set}{\mathbf{Set}}
\newcommand{\Art}{\mathbf{Art}}
\newcommand{\Alg}{\mathbf{Alg}}
\newcommand{\DGArt}{\mathbf{DGArt}}
\newcommand{\CDGA}{\mathbf{CDGA}}
\newcommand{\CGA}{\mathbf{CGA}}
\newcommand{\DGLA}{\mathbf{DGLA}}
\newcommand{\solose}{\Rightarrow}
\newcommand{\PSI}{\Psi\mathbf{Sch}_I(\mathbf{M})}
\newcommand{\PSJ}{\Psi\mathbf{Sch}_J(\mathbf{M})}
\newcommand{\Sym}{\mbox{Sym}}

\title{A DG-Enhancement of $\D(\QCoh(X))$ with applications in deformation theory}
\date{August 15, 2018}
\author{Francesco Meazzini}
\address{\newline
Universit\`a degli studi di Roma La Sapienza,\hfill\newline
Dipartimento di Matematica \lq\lq Guido
Castelnuovo\rq\rq,\hfill\newline
P.le Aldo Moro 5,
I-00185 Roma, Italy.}
\email{meazzini@mat.uniroma1.it}


\maketitle

\begin{abstract}
It is well-known that DG-enhancements of $\D(\QCoh(X))$ are all equivalent to each other, see~\cite{LO}. Here we present an explicit model which leads to applications in deformation theory. In particular, we shall describe three models for derived endomorphisms of a quasi-coherent sheaf $\sF$ on a finite-dimensional Noetherian separated scheme (even if $\sF$ does not admit a locally free resolution). Moreover, these complexes are endowed with DG-Lie algebra structures, which we prove to control infinitesimal deformations of $\sF$.
\end{abstract}

\tableofcontents

\section{Introduction}


A classical problem in deformation theory concerns the study of infinitesimal deformations of a quasi-coherent sheaf $\sF$ on a scheme $X$ over a field $\K$ of characteristic $0$.
Deformations up to isomorphisms define a functor $\Def_{\sF}\colon \Art_{\K}\to\Set$, where $\Art_{\K}$ denotes the category of local Artin $\K$-algebras with residue field $\K$. The classical approach is based on a (finite) locally free resolution $\sE\to\sF$, which for instance exists if $X$ is smooth projective. In fact, a deformation of $\sF$ can be understood as the data of local deformations of $\sE$ together with suitable gluing conditions.
It is proven in~\cite{FIM} that $\Def_{\sF}$ is controlled by the DG-Lie algebra of global sections of an acyclic resolution of the sheaf $\sE nd^{\ast}(\sE)$ in the sense of~\cite{GM,Man}. In particular, it is well-known that $T^1\Def_{\sF}\cong\Ext^1(\sF,\sF)$ and obstructions are contained in $\Ext^2(\sF,\sF)$. This highlights the considerable role of derived endomorphisms $\R\End(\sF)$, and the importance of being able to compute its cohomology $\Ext^{\ast}(\sF,\sF)$. Classically, $\R\End(\sF)$ is defined (up to quasi-isomorphisms) as the complex $\Hom^{\ast}_{\Oh_X}(\sF,\sI)$ for any injective resolution $\sF\to\sI$. Unfortunately, despite the outstanding fact that injective resolutions always exist, it is often very hard to describe them.

Here comes the aim of this paper to present another approach to compute $\R\End(\sF)$ when dealing with concrete geometric situations, always trying to keep the exposition as clear as possible with the attempt to reduce the use of simplicial and model category techniques at minimum.

The main tool is the introduction of the category $\Mod(A_{\pal})$ of modules over (the diagram $A_{\pal}$ representing) a separated $\K$-scheme $X$. Fix an open affine covering $\sU=\{U_h\}$ for $X$, then the associated diagram $A_{\pal}$ with respect to $\sU$ is defined as
\[ A_{\pal}\colon \sN\to \Alg_{\K} \; ,\qquad \qquad \alpha\mapsto A_{\alpha}=\Gamma(U_{\alpha},\Oh_X) \]
where $\sN=\{ \alpha=\{h_0,\dots,h_k\}\,\vert\,U_{\alpha}=U_{h_0}\cap\dots U_{h_k}\neq\emptyset \}$ is the nerve of $\sU$.
An $A_{\pal}$-module $\sG$ can be understood as the following data 
\begin{enumerate}
\item a DG-module $\sG_{\alpha}$ over $A_{\alpha}$ for every $\alpha$ in the nerve $\sN$ of $\sU$,
\item a morphism $g_{\alpha\beta}\colon \sG_{\alpha}\otimes_{A_{\alpha}}A_{\beta}\to \sG_{\beta}$ of $A_{\beta}$-modules, for every $\alpha\subseteq\beta$ in $\sN$,
\end{enumerate}
satisfying the \emph{cocycle condition}, see Definition~\ref{def.pseudo-module}.
Similar notions were considered in~\cite{EE,FK,GS,Sto}. Taking advantage of the standard projective model structure on DG-modules, the category $\Mod(A_{\pal})$ will be endowed with a (cofibrantly generated) model structure, see Theorem~\ref{thm.modelpseudomodules}, where weak equivalences are pointwise quasi-isomorphisms. The above model structure can be seen as a geometric example of an abstract recent result obtained in~\cite{Bal}. In order to work with quasi-coherent sheaves, we need a (homotopical) version of quasi-coherence for $A_{\pal}$-modules: $\sG$ is called quasi-coherent if all the maps $g_{\alpha\beta}$ introduced above are quasi-isomorphisms, see Definition~\ref{def.qcohpseudo-module}.
To the author knowledge the last definition does not appear in the existing literature, a part for the case of non-graded modules for which the theory is carried out in~\cite{EE,Sto}. Now, denote by $\Ho(\QCoh(A_{\pal}))$  the category of quasi-coherent $A_{\pal}$-modules localized with respect to the weak equivalences: Theorem~\ref{thm.equivalence} states that this is equivalent to the (unbounded) derived category of quasi-coherent sheaves on $X$, hence leading to an explicit description of a DG-enhancement of $\D(\QCoh(X))$, see Corollary~\ref{corollary.enhancement}. It is worth to notice that some of the functors involved in Section~\ref{section.derivedlowershriek} have been somehow already considered in the literature, see~\cite{Hin04,HS}. Moreover a result similar to the equivalence of Theorem~\ref{thm.equivalence} was partially proven in~\cite[Proposition 2.28]{BF}.

In~\cite{LO} it was shown the uniqueness of DG-enhancements for the derived category of a suitable Grothendieck category up to equivalence. In particular, this applies to $\D(\QCoh(X))$ under some mild hypothesis on $X$ (e.g. if $X$ is a quasi-projective $\K$-scheme); therefore the mere description of a DG-enhancement is not very exciting.

On the other hand, our construction turns out to be very useful when dealing with derived endomorphisms of a quasi-coherent sheaf $\sF$ of $\Oh_X$-modules.
In fact, the category of $A_{\pal}$-modules allows to easily describe $\R\End(\sF)$ in terms of a cofibrant replacement of $\sF$, see Theorem~\ref{thm.REndpseudomodules}. Moreover, Example~\ref{example.replacementLocallyFree} shows the feasibility of the computation of such cofibrant replacement in interesting cases.
In Section~\ref{section.REnd} we propose two more models for $\R\End(\sF)$: the first is again in terms of a cofibrant replacement in the model category of $A_{\pal}$-modules and involves the Thom-Whitney totalization, Corollary~\ref{corollary.diagram}, while the other assumes the existence of a locally free resolution for $\sF$, Theorem~\ref{thm.locallyfreeREnd}.

The last section is devoted to our main application in deformation theory; in particular, we deal with the functor $\Def_{\sF}\colon\Art_{\K}\to\Set$ of classical infinitesimal deformations of $\sF$.
Recall that since the eighties the leading principle in deformation theory (due to Quillen, Deligne, Drinfeld, Kontsevich...) states that any deformation problem is controlled by a DG-Lie algebra via Maurer-Cartan solutions modulo gauge equivalence, see~\cite{GM,Man,Manin17}.
Around 2010 this was formally proven independently by Lurie~\cite[Theorem 5.3]{Lur} and Pridham~\cite[Theorem 4.55]{Pri}; it is dutiful to mention that partial results in this direction where previously obtained by Hinich and Manetti, see~\cite{Hin01,Man2,Pri} and references therein.
In Section~\ref{section.deformations} we adopt this point of view proving that the three complexes representing $\R\End(\sF)$ described in Section~\ref{section.REnd} are all equipped with a DG-Lie algebra structure, and each of them controls $\Def_{\sF}$ via Maurer-Cartan elements modulo gauge equivalence. In particular, we give two proofs of this fact: the first (Section~\ref{section.FIM}) involves the semicosimplicial machinery together with standard arguments of descent of the Deligne groupoid, while the second (Section~\ref{section.deformationpseudomod}) relies on a direct computation in $\Mod(A_{\pal})$.

A remarkable fact is that our descriptions of $\R\End(\sF)$ in terms of $A_{\pal}$-modules does not require the existence of a locally free resolution for $\sF$, since cofibrant replacements always exist. Hence we recover that $T^1\Def_{\sF}\cong\Ext^1(\sF,\sF)$ and that obstructions are contained in $\Ext^2(\sF,\sF)$ only assuming $X$ to be a finite-dimensional Noetherian separated $\K$-scheme.

\bigskip

\textbf{Acknowledgements.} I am deeply in debt with Marco Manetti, who encouraged me to draw up this paper and supported me with several useful suggestions.

\bigskip

\section{Preliminaries and notation}

This short introductory section aims to fix the geometric framework where we shall work throughout all the paper, and to briefly recall some basic constructions.

We work on a fixed finite-dimensional Noetherian separated scheme $X$ over a field $\K$ of characteristic $0$. Actually, the assumption on the characteristic of $\K$ will be necessary only in Section~\ref{section.REnd} and Section~\ref{section.deformations}, where applications to algebraic geometry will be discussed. For simplicity of exposition we shall work over $\K$ throughout all the paper, although the results of the first sections hold for schemes over $\Z$. Moreover, we fix an open affine covering $\sU=\{U_h\}_{h\in H}$ together with its \emph{nerve}
\[ \sN = \left\{ \{h_0,\dots,h_k\} \,\vert\, U_{h_0}\cap\dots\cap U_{h_k}\neq\emptyset \right\} \]
which carries a \emph{degree function} $\deg\colon\sN\to\N$ defined by $\deg(\{h_0,\dots,h_k\}) = k$. Moreover, for every $\alpha=\{h_0,\dots, h_k\}\in\sN$ we denote by $U_{\alpha}$ the intersection $U_{h_0}\cap\dots\cap U_{h_k}$, and by $A_{\alpha}=\Gamma(U_{\alpha},\Oh_X)$.
Recall that since $X$ is assumed to be separated, then $U_{\alpha}$ is in fact affine for every $\alpha\in\sN$.
The nerve $\sN$ is a partially ordered set where $\alpha\leq\beta$ if and only if $\alpha\subseteq\beta$; notice that if $\alpha\leq\beta$ then $U_{\beta}\subseteq U_{\alpha}$ so that there exists a flat map of $\K$-algebras $A_{\alpha}\to A_{\beta}$ satisfying $A_{\beta}\cong A_{\beta}\otimes_{A_{\alpha}}A_{\beta}$.
Hence, once we have fixed $\sU$, the scheme $X$ can be represented by the diagram
\[ A_{\pal}\colon \sN\to\Alg_{\K} \; , \qquad \qquad \alpha\mapsto A_{\alpha} \]
where $A_{\alpha}\to A_{\beta}$ is the opposite map of $\Spec(A_{\beta})\to\Spec(A_{\alpha})$ for every $\alpha\leq\beta$ in $\sN$.

For any open subset $U\subseteq X$ let $\DGMod(\Oh_U)$ be the category of (unbounded) complexes of $\mathcal{O}_U$-modules, and by $\QCoh(U)$ the full subcategory of complexes of quasi-coherent sheaves.

For every inclusion $i\colon U\to V$ between open subsets there are three associated functors:
\[ i_{!},i_{\ast}\colon \DGMod(\mathcal{O}_U)\to \DGMod(\mathcal{O}_V),\qquad 
i^{\ast}\colon \DGMod(\mathcal{O}_V)\to \DGMod(\mathcal{O}_U)\,.\]
Recall that $i^{\ast}\sG=\sG\vert_{U}$ because $\Oh_V\vert_U=\Oh_U$, and 
$i_!\sF$ is the sheaf associated to the presheaf $i(\sF)$ defined by
\[ \begin{cases}
i(\sF)(W)=\sF(W) & \text{ if }W\subseteq U \\
i(\sF)(W)=0 & \text{ otherwise}.
\end{cases} \]

The obvious retraction $i(\sF)\to i_{\ast}(\sF)\to i(\sF)$ of presheaves gives a retraction of sheaves 
$i_!\sF\to i_{\ast}\sF\to i_!\sF$ and then a retraction of functors $i_!\to i_{\ast}\xrightarrow{r}i_!$.
Notice also that for every $\sG\in \DGMod(\mathcal{O}_V)$ there exists an injective morphism 
\[ i_!i^{\ast}\sG\to \sG\]
and therefore a natural morphism given by composition with the retraction $r$
\[  i_{\ast}i^{\ast}\sG\to \sG \; , \]
which is an isomorphism on stalks over every $x\in U$, and $0$ over $x\notin U$.

If $\sF$ and $\sG$ are complexes of quasi-coherent sheaves, then also $i_{\ast}\sF$ and $i^{\ast}\sG$ are so, see e.g.~\cite[Proposition 5.8]{Har}. This is not true in general for $i_!\sF$, see e.g.~\cite[Example 5.2.3]{Har}.

Recall that in the above notation, if $U$ is affine then the functor $i_{\ast}\colon \QCoh(U)\to \QCoh(V)$ is exact.

\bigskip

\section{The model category of $A_{\pal}$-modules}

The aim of this section is to introduce the category $\Mod(A_{\pal})$ of modules over the diagram representing the scheme $X$ with respect to a fixed open affine covering. Moreover, we shall endow $\Mod(A_{\pal})$ with a model structure, see Theorem~\ref{thm.modelpseudomodules}, which will lead us to our main result in Section~\ref{section.equivalence}.

In the following, for every ring $R$ we denote by $\DGMod(R)$ the category of DG-modules over $R$. Recall that the diagram $A_{\pal}\colon \sN\to \Alg_{\K}$ is defined by $A_{\alpha}=\Gamma(U_{\alpha},\Oh_X)$.

\begin{definition}\label{def.pseudo-module}
An $A_{\pal}$\textbf{-module} $\mathcal{F}$ over $X$, with respect to the fixed covering $\sU$, consists of the following data:
\begin{enumerate}
\item an object $\sF_{\alpha}\in\DGMod(A_{\alpha})$, for every $\alpha\in\sN$,
\item a morphism $f_{\alpha\beta}\colon \sF_{\alpha}\otimes_{A_{\alpha}}A_{\beta}\to \sF_{\beta}$ in $\DGMod(A_{\beta})$, for every $\alpha\leq\beta$ in $\sN$,
\end{enumerate}
satisfying the \emph{cocycle condition} $f_{\beta\gamma}\circ\left(f_{\alpha\beta}\otimes_{A_{\beta}}A_{\gamma}\right) = f_{\alpha\gamma}$, for every $\alpha\leq\beta\leq\gamma$ in $\sN$.
\end{definition}

In the setting of Definition~\ref{def.pseudo-module}, the map $f_{\alpha\beta}\colon \sF_{\alpha}\otimes_{A_{\alpha}}A_{\beta}\to \sF_{\beta}$ in $\DGMod(A_{\beta})$ is equivalent to its adjoint morphism $\sF_{\alpha}\to\sF_{\beta}$ in $\DGMod(A_{\alpha})$, where the $A_{\alpha}$-module structure on $\sF_{\beta}$ is induced via $A_{\alpha}\to A_{\beta}$.

For instance, to any sheaf $\sG$ of $\Oh_X$-modules it is associated the $A_{\pal}$-module $\Upsilon^{\ast}\sG$ defined as
\[ \left(\Upsilon^{\ast}\sG\right)_{\alpha} = \sG(U_{\alpha})\in\DGMod(A_{\alpha}) \qquad \text{ and } \qquad g_{\alpha\beta}\colon\sG(U_{\alpha})\otimes_{A_{\alpha}}A_{\beta}\to \sG(U_{\beta}) \]
for every $\alpha\leq\beta$ in $\sN$, where the map $g_{\alpha\beta}$ is induced by the restriction map of the sheaf $\sG$.

\begin{definition}\label{definition.pseudomap}
A \textbf{morphism of $A_{\pal}$-modules} $\varphi\colon \mathcal{F}\to \mathcal{G}$ over $X$ consists of the following data:
\begin{enumerate}
\item a morphism $\varphi_{\alpha}\colon \sF_{\alpha}\to \sG_{\alpha}$ in $\DGMod(A_{\alpha})$, for every $\alpha\in \sN$,
\item for every $\alpha\leq\beta$ in $\sN$, the diagram
\[ \xymatrix{	\sF_{\alpha}\otimes_{A_{\alpha}}A_{\beta} \ar@{->}[r]^{\varphi_{\alpha}} \ar@{->}[d]_{f_{\alpha\beta}} & \sG_{\alpha}\otimes_{A_{\alpha}}A_{\beta} \ar@{->}[d]^{g_{\alpha\beta}} \\
\sF_{\beta} \ar@{->}[r]_{\varphi_{\beta}} & \sG_{\beta}	} \]
commutes in $\DGMod(A_{\beta})$.
\end{enumerate}
The set of morphisms between $\sF$ and $\sG$ is denoted by $\Hom_{A_{\pal}}(\sF,\sG)$.
\end{definition}

Recall that for any ring $R$ and any pair of DG-modules $M,N\in\DGMod(R)$ it is defined total-Hom complex $\Hom_R^{\ast}(M,N)$ as follows:
\[ \Hom_R^{p}(M,N)=\prod\limits_{n\in\Z}\Hom_R(M^n,N^{n+p}) \; ,  \qquad \partial_{\Hom}^p\colon (f^n)_{n\in\Z}\mapsto (f^{n+1}d_N^n-(-1)^{p} d_N^{n+p}f^n)_{n\in\Z} \; .\]

\begin{definition}\label{definition.DGpseudomap}
\textbf{$\ast$-morphisms} between $A_{\pal}$-modules $\mathcal{F}$ and $\mathcal{G}$ over $X$ are defined by:
\[ \Hom_{A_{\pal}}^{\ast}(\sF,\sG)\subseteq\prod_{\alpha\in\sN}\Hom_{A_{\alpha}}^{\ast}(\sF_{\alpha},\sG_{\alpha}) \]
where $\{\varphi_{\alpha}\}_{\alpha\in\sN}$ belongs to $\Hom_{A_{\pal}}^{\ast}(\sF,\sG)$ if the diagram
\[ \xymatrix{	\sF_{\alpha}\otimes_{A_{\alpha}}A_{\beta} \ar@{->}[r]^{\varphi_{\alpha}} \ar@{->}[d]_{f_{\alpha\beta}} & \sG_{\alpha}\otimes_{A_{\alpha}}A_{\beta} \ar@{->}[d]^{g_{\alpha\beta}} \\
\sF_{\beta} \ar@{->}[r]_{\varphi_{\beta}} & \sG_{\beta}	} \]
commutes for every $\alpha\leq\beta\in\sN$.
\end{definition}

Notice that $\Hom_{A_{\pal}}(\sF,\sG)$ are precisely the $0$-cocycles of the complex $\Hom_{A_{\pal}}^{\ast}(\sF,\sG)$, whose differential is the inherited (graded) commutator.
We shall denote by $\Mod(A_{\pal})$ the category of $A_{\pal}$-modules, with morphisms of $A_{\pal}$-modules as in Definition~\ref{definition.pseudomap}.
We shall denote by $\DGPM(A_{\pal})$ the DG-category of $A_{\pal}$-modules, with $\ast$-morphisms as in Definition~\ref{definition.DGpseudomap}. Since the covering $\sU$ is assumed to be fixed at the beginning, we do not emphasise the dependence on it.

Recall that by~\cite{DS,Hov99,QuillenHA} for every $\alpha\in\sN$ the category $\DGMod(A_{\alpha})$ is endowed with a model structure where
\begin{itemize}
\item weak equivalences are quasi-isomorphisms,
\item fibrations are degreewise surjective morphisms,
\item every object is fibrant
\item $C\in\DGMod(A_{\alpha})$ is cofibrant if and only if for every cospan $C\xrightarrow{f} D\xleftarrow{g} E$ with $g$ a surjective quasi-isomorphism there exists a lifting $h\colon C\to E$.
\item cofibrations are degreewise split injective morphisms with cofibrant cokernel.
\end{itemize}
Recall that if a complex $Y\in\DGMod(A_{\alpha})$ is bounded above then it is cofibrant if and only if it is degreewise projective, see~\cite[Lemma 2.3.6]{Hov99}. 

Our next goal is to endow the category $\Mod(A_{\pal})$ with a model structure. To this aim, we first need preliminary results. Fix $\sF\in\PM(A_{\pal})$ and $\alpha\in\sN$; define the \textbf{latching module} of $\sF$ at $\alpha$ to be
\[ L_{\alpha}\sF=\colim_{\gamma<\alpha}\left(\sF_{\gamma}\otimes_{A_{\gamma}}A_{\alpha}\right) \]
and notice that there exists a natural map $L_{\alpha}\sF\to \sF_{\alpha}$. This allows us to define the full subcategory of cofibrant $A_{\pal}$-modules: an $A_{\pal}$-module $\sF\in\PM(A_{\pal})$ is called \textbf{cofibrant} if for every $\alpha\in\sN$ the latching map $L_{\alpha}\sF\to \sF_{\alpha}$ is a cofibration in $\DGMod(A_{\alpha})$.
Cofibrant $A_{\pal}$-modules define full subcategories $\PM(A_{\pal})^c\subseteq\PM(A_{\pal})$ and $\PM^{\ast}(A_{\pal})^c\subseteq\PM^{\ast}(A_{\pal})$; in particular, we shall prove that the category $\PM^{\ast}(A_{\pal})^c$ is a DG-enhancement (in the sense of~\cite{LO}) of the unbounded derived category $\D(\QCoh(X))$, see Corollary~\ref{corollary.enhancement}. 

\begin{remark}\label{rmk.simplicialhomology}
Let $\{U_h\}_{h\in H}$ be an open cover of $X$ and let $\sN$ be its nerve. Choose a total order on $H$; then to every $\alpha\in\sN$ it is associated the abstract oriented simplicial complex $\mathcal{P}(\alpha)$, whose faces are the subsets of $\alpha$. Moreover, denote by $C_{\alpha}$ the corresponding chain complex. Recall that $C_{\alpha}$ in degree $r$ is the free abelian group of rank $\binom{\deg(\alpha)+1}{r+1}$, and its homology is given by: $H^0(C_{\alpha})=\Z$ and $H^j\left(C_{\alpha}\right)=0$ for every $j\neq0$. Now consider the category $\Ch(\Z)$ of chain complexes of abelian groups; we define the diagram
\[ C\colon \sN\to \Ch(\Z); \qquad \qquad \alpha\mapsto C_{\alpha} \]
where for every $\alpha\leq\beta$ in $\sN$ the map $C_{\alpha}\to C_{\beta}$ is the natural inclusion. Notice that there is a short exact sequence
\[ 0 \to \colim_{\gamma<\alpha}C_{\gamma} \xrightarrow{\iota_{\alpha}} C_{\alpha} \to \coker(\iota_{\alpha}) \to 0 \]
where $\coker(\iota_{\alpha})^{\deg(\alpha)}=\Z$ and $\coker(\iota_{\alpha})^j=0$ for every $j\neq \deg(\alpha)$.
\end{remark}

\begin{example}[Cofibrant $A_{\pal}$-module associated to $\Oh_X$]\label{example.cofibrantOX}
To the scheme $X$ it is associated the $A_{\pal}$-module $\sQ_X\in\PM(A_{\pal})$ defined as
\[ \sQ_{X,\alpha}^r = C_{\alpha}^{-r}\otimes_{\Z}A_{\alpha} \qquad \text{ and } \qquad d_{\sQ_X}^r = d_{C_{\alpha}}^{-r} \otimes \Id_{A_{\alpha}} \]
for every $r\in\Z$ and every $\alpha\in\sN$. Therefore $\sQ_{X,\alpha}\in\DGMod(A_{\alpha})$, and for every $\alpha\leq\beta$ the map
\[ \sQ_{X,\alpha}\otimes_{A_{\alpha}}A_{\beta}\to\sQ_{X,\beta} \]
is induced by the natural inclusion $C_{\alpha}\to C_{\beta}$. For simplicity of notation we shall denote by $C_{\alpha}^{op}$ the cochain complex defined by $(C_{\alpha}^{op})^r=C_{\alpha}^{-r}$ and $d_{C_{\alpha}^{op}}^r=d_{C_{\alpha}}^{-r}$, $r\in\Z$; hence $\sQ_{X,\alpha}=C_{\alpha}^{op}\otimes_{\Z}A_{\alpha}$ for every $\alpha\in\sN$. Notice that by Remark~\ref{rmk.simplicialhomology} for every $\alpha\in\sN$ we have a short exact sequence
\[ 0 \to L_{\alpha}\sQ_X \xrightarrow{\iota_{\alpha}\otimes \Id_{A_{\alpha}}} \sQ_{X,\alpha} \to \coker(\iota_{\alpha})\otimes_{\Z}A_{\alpha} \to 0\]
so that the latching map $\iota_{\alpha}\otimes \Id_{A_{\alpha}}$ is degreewise injective and its cokernel is zero except for degree $\deg(\alpha)$. Finally, since $\coker(\iota_{\alpha}\otimes \Id_{A_{\alpha}})^{\deg(\alpha)} = A_{\alpha}$ is a free (hence projective) $A_{\alpha}$-module, then the latching map is in fact a cofibration in $\DGMod(A_{\alpha})$ by~\cite[Lemma 2.3.6]{Hov99}. This proves that $\sQ_X$ is a cofibrant $A_{\pal}$-module.
\end{example}

Cofibrant $A_{\pal}$-modules play an important role also in the application to deformation theory that we will describe in Section~\ref{section.deformations}.
With this in mind, the crucial notion is the following: a \textbf{cofibrant replacement} for a given $A_{\pal}$-module $\sF\in\PM(A_{\pal})$ is a morphism $\sQ\to\sF$ in $\PM(A_{\pal})$ such that
\begin{description}
\item[1] $\sQ$ is a cofibrant $A_{\pal}$-module,
\item[2] the map $\sQ_{\alpha}\to\sF_{\alpha}$ is a surjective quasi-isomorphism for every $\alpha\in\sN$.
\end{description}
Cofibrant replacements are usually not unique and should be though as global resolutions of the $A_{\pal}$-module $\sF$. 

\begin{example}[Cofibrant replacement for the structure sheaf $\Oh_X$]\label{example.replacementOX}
As already noticed, to any sheaf $\sG$ of $\Oh_X$-modules it is associated an $A_{\pal}$-module $\Upsilon^{\ast}\sG\in\PM(A_{\pal})$. In particular, $\Upsilon^{\ast}\Oh_X\in\PM(A_{\pal})$ is defined as $\left(\Upsilon^{\ast}\Oh_X\right)_{\alpha} = A_{\alpha}$ on every $\alpha\in\sN$, and the map $\left(\Upsilon^{\ast}\Oh_X\right)_{\alpha}\otimes_{A_{\alpha}}A_{\beta}\to \left(\Upsilon^{\ast}\Oh_X\right)_{\beta}$ is the identity for every $\alpha\leq\beta$.

In Example~\ref{example.cofibrantOX} we constructed the cofibrant $A_{\pal}$-module $\sQ_X\in\PM(A_{\pal})$. Notice that by Remark~\ref{rmk.simplicialhomology} the set of maps $\{C_{\alpha}\to H^0(C_{\alpha}) = \Z\}_{\alpha\in\sN}$ induce a morphism $\sQ_X\to\Upsilon^{\ast}\Oh_X$ in $\PM(A_{\pal})$ which is a cofibrant replacement. In fact, by the flatness of the map $A_{\alpha}\to A_{\beta}$ it follows that
\[ \pi_{\alpha}\colon \sQ_{X,\alpha} = C_{\alpha}^{op} \otimes_{\Z} A_{\alpha} \to A_{\alpha} = \left( \Upsilon^{\ast}\Oh_X \right)_{\alpha} \]
is a surjective quasi-isomorphism for every $\alpha\in\sN$.
\end{example}

\begin{example}[Cofibrant replacement for a locally free sheaf]\label{example.replacementLocallyFree}
Consider a locally free sheaf $\sE$ on $X$, and take a cover $\{U_h\}_{h\in H}$ such that $\sE\vert_{U_{\alpha}}$ is a free $A_{\alpha}$-module for every $\alpha\in\sN$.
Since for every $\alpha\in\sN$ the (DG-)module $\Upsilon^{\ast}\sE_{\alpha}=\sE(U_{\alpha})$ is concentrated in degree $0$, it is cofibrant in $\DGMod(A_{\alpha})$ by~\cite[Lemma 2.3.6]{Hov99}. Nevertheless, the latching maps need not to be cofibrations in general; hence $\Upsilon^{\ast}\sE$ is an example of an $A_{\pal}$-module which is pointwise cofibrant but not (globally) cofibrant. Following Example~\ref{example.replacementOX} we can explicitly construct a cofibrant replacement $\sQ_{\sE}\to\Upsilon^{\ast}\sE$:
\begin{itemize}
\item $\sQ_{\sE,\alpha} = \sQ_{X,\alpha} \otimes_{A_{\alpha}} \sE(U_{\alpha})$ for every $\alpha\in\sN$,
\item for every $\alpha\leq\beta$ in $\sN$ the morphism $\sQ_{\sE,\alpha}\otimes_{A_{\alpha}}A_{\beta}\to \sQ_{\sE,\beta}$ is induced by the corresponding restriction map of $\sE$,
\item the morphism $\sQ_{\sE,\alpha}\to \sE(U_{\alpha}) = (\Upsilon^{\ast}\sE)_{\alpha}$ is defined as $\pi_{\alpha}\otimes\Id_{\sE(U_{\alpha})}$ for every $\alpha\in\sN$.
\end{itemize}
By Example~\ref{example.replacementOX} $\pi\colon \sQ_X\to\Upsilon^{\ast}\Oh_X$ is a cofibrant replacement; therefore the map $\pi\otimes\Id \colon\sQ_{\sE}\to\Upsilon^{\ast}\sE$ is a cofibrant replacement for $\Upsilon^{\ast}\sE$.
\end{example}

Now fix $\alpha\in\sN$; define $\mathbf{R}_{\alpha}=\{\gamma\in\sN\,\vert\,\gamma<\alpha\}$ and recall that the category of diagrams $\DGMod(A_{\alpha})^{\mathbf{R}_{\alpha}}$ is endowed with the Reedy model structure where Reedy weak equivalences and Reedy fibrations are detected pointwise. Namely, a map $f\colon Y\to Z$ in $\DGMod(A_{\alpha})^{\mathbf{R}_{\alpha}}$ is a Reedy weak equivalence (respectively: Reedy fibration) if and only if $f_{\gamma}\colon Y_{\gamma}\to Z_{\gamma}$ is a quasi-isomorphism (respectively: degreewise surjective) for every $\gamma<\alpha$. Moreover, $f$ is a Reedy cofibration if and only if the map
\[  \colim_{\beta<\gamma}Z_{\beta} \coprod\limits_{\colim\limits_{\beta<\gamma}Y_{\beta}} Y_{\gamma} \to Z_{\gamma} \]
is a cofibration in $\DGMod(A_{\gamma})$ for every $\gamma\in\mathbf{R}_{\alpha}$, see~\cite[Theorem 16.3.4]{Hir03}.

We have functors $\res_{\alpha}\colon \PM(A_{\pal}) \to \DGMod(A_{\alpha})^{\mathbf{R}_{\alpha}}$ defined by
\[ (\res_{\alpha}\sF)_{\gamma} = \sF_{\gamma}\otimes_{A_{\gamma}}A_{\alpha}\; , \qquad \gamma<\alpha \]
for every $\sF\in\PM(A_{\pal})$.

\begin{lemma}\label{lem.trivialcofibration}
For every morphism $\varphi\colon\sF\to \sG$ in $\PM(A_{\pal})$ the following conditions are equivalent.
\begin{enumerate}
\item For every $\alpha\in\sN$, the morphism $\varphi_{\alpha}\colon \sF_{\alpha}\to \sG_{\alpha}$ is a quasi-isomorphism in $\DGMod(A_{\alpha})$, and the natural morphism
\[ L_{\alpha}\sG\amalg_{(L_{\alpha}\sF)} \sF_{\alpha} \longrightarrow \sG_{\alpha} \]
is a cofibration in $\DGMod(A_{\alpha})$.
\item For every $\alpha\in\sN$, the natural morphism
\[ L_{\alpha}\sG\amalg_{(L_{\alpha}\sF)} \sF_{\alpha} \longrightarrow \sG_{\alpha} \]
is a trivial cofibration in $\DGMod(A_{\alpha})$.
\end{enumerate}
\begin{proof}
Fix $\alpha\in\sN$ and consider the following diagram
\[ \xymatrix{	L_{\alpha}\sF \ar@{->}[r] \ar@{->}[d] & \sF_{\alpha} \ar@{->}[d] \ar@{->}@/^15pt/[rdd]^{\varphi_{\alpha}} & \\
L_{\alpha}\sG \ar@{->}[r] \ar@{->}@/_15pt/[rrd] & L_{\alpha}\sG\amalg_{(L_{\alpha}\sF)} \sF_{\alpha} \ar@{.>}[dr]|-{\psi} & \\
& & \sG_{\alpha}	\, .} \]

Now define two diagrams in $\DGMod(A_{\alpha})^{\mathbf{R}_{\alpha}}$ as $Y=\res_{\alpha}\sF \mbox{ and } Z=\res_{\alpha}\sG$, and notice that if either $(1)$ or $(2)$ holds the morphism $Z\to Y$ induced by $\varphi$ is a Reedy cofibration, since colimits commute with coproducts. Moreover, by~\cite[Theorem 15.3.15]{Hir03} it follows that $Y\to Z$ is a Reedy weak equivalence if either $(1)$ or $(2)$ holds, so that the vertical morphisms in the diagram above are trivial cofibrations in $\DGMod(A_{\alpha})$; in fact $\colim\colon \DGMod(A_{\alpha})^{\mathbf{R}_{\alpha}} \to \DGMod(A_{\alpha})$ is a left Quillen functor and trivial cofibrations are closed under pushouts. Therefore, $\varphi_{\alpha}$ is a weak equivalence if and only if $\psi$ is so, because of the \emph{2 out of 3} axiom.
\end{proof}
\end{lemma}

The following result endows the category $\PM(A_{\pal})$ with a model structure where the class of cofibrant objects coincides with cofibrant $A_{\pal}$-modules defined above.

\begin{theorem}[\textbf{Model structure on $A_{\pal}$-modules}]\label{thm.modelpseudomodules}
The category of $A_{\pal}$-modules over $X$ is endowed with a model structure, where a morphism $\sF\to \sG$ in $\PM(A_{\pal})$ is
\begin{enumerate}
\item a weak equivalence if and only if the morphism $\sF_{\alpha}\to \sG_{\alpha}$ is a quasi-isomorphism in $\DGMod(A_{\alpha})$ for every $\alpha\in \sN$,
\item a fibration if and only if the morphism $\sF_{\alpha}\to \sG_{\alpha}$ is surjective in $\DGMod(A_{\alpha})$ for every $\alpha\in \sN$,
\item a cofibration if and only if the natural morphism
\[ L_{\alpha}\sG\amalg_{(L_{\alpha}\sF)} \sF_{\alpha} \longrightarrow \sG_{\alpha} \]
is a cofibration in $\DGMod(A_{\alpha})$ for every $\alpha\in\sN$.
\end{enumerate}
\begin{proof}
It is sufficient to prove that $\PM(A_{\pal})$ with the classes defined in the statement satisfies the axioms of a model category. First notice that the category $\PM(A_{\pal})$ is complete and cocomplete since limits and colimits are taken pointwise. Moreover, the class of weak equivalences satisfies the \emph{2 out of 3} axiom by definition.

The closure with respect to retracts holds since if $\sF\to\sG$ is a retract of $\sF'\to \sG'$ in the category of maps of $\PM(A_{\pal})$, then the natural morphism $L_{\alpha}\sG\amalg_{(L_{\alpha}\sF)} \sF_{\alpha} \longrightarrow \sG_{\alpha}$ is a retract of the natural morphism $L_{\alpha}\sG'\amalg_{(L_{\alpha}\sF')} \sF'_{\alpha} \longrightarrow \sG'_{\alpha}$ in the category of maps of $\DGMod(A_{\alpha})$, for every $\alpha\in\sN$.

In order to show that the \emph{lifting} axiom holds, observe that a morphism $\mathcal{F}\to\mathcal{G}$ is a trivial cofibration in $\PM(A_{\pal})$ if and only if for every $\alpha\in\sN$ the natural morphism $L_{\alpha}\sG\amalg_{(L_{\alpha}\sF)} \sF_{\alpha} \longrightarrow \sG_{\alpha}$ is a trivial cofibration in $\DGMod(A_{\alpha})$, see Lemma~\ref{lem.trivialcofibration}. Therefore the required lifting can be constructed inductively on the degree of $\alpha$.

The \emph{factorization} axiom can be proved inductively as follows. Take a morphism $\varphi\colon\mathcal{F}\to \mathcal{G}$, we need to define (functorial) factorizations $\mathcal{F}\to \mathcal{Q}\to\mathcal{G}$ in $\PM(A_{\pal})$ as a cofibration (respectively, trivial cofibration) followed by a trivial fibration (respectively, fibration). Now, fix $\alpha\in\sN$ of degree $d$ and suppose $\varphi_{\gamma}$ has been factored for all $\gamma\in\sN$ of degree less that $d$. Consider a (functorial) factorization of the natural morphism
\[ L_{\alpha}\sG\amalg_{(L_{\alpha}\sF)} \sF_{\alpha}  \longrightarrow \sQ_{\alpha}\longrightarrow \sG_{\alpha} \]
in $\DGMod(A_{\alpha})$ as a cofibration (respectively, trivial cofibration) followed by a trivial fibration (respectively, fibration). Lemma~\ref{lem.trivialcofibration} implies that $\mathcal{Q}$ satisfies the required properties by construction.
\end{proof}
\end{theorem}

\begin{remark}
A morphism $f\colon\sF\to \sG$ in $\PM(A_{\pal})$ is a weak equivalence (respectively: fibration, cofibration) with respect to the model structure of Theorem~\ref{thm.modelpseudomodules} if and only if for every $\alpha\in\sN$ the induced morphism $\res_{\alpha}(f)$ is a Reedy weak equivalence (respectively: Reedy fibration, Reedy cofibration) in $\DGMod(A_{\alpha})^{\mathbf{R}_{\alpha}}$. The claim immediately follows by the flatness of the map $A_{\beta}\to A_{\gamma}$ for every $\beta\leq\gamma$.
\end{remark}

\begin{remark}\label{rmk.boundedpseudomodule}
For any $\alpha\in\sN$, consider the full subcategory $\DGMod^{\le 0}(A_{\alpha})\subseteq\DGMod(A_{\alpha})$ whose objects are complexes concentrated in non-positive degrees.
This is endowed with a model structure where
\begin{itemize}
\item weak equivalences are quasi-isomorphisms,
\item fibrations are surjections in negative degrees,
\item cofibrations are degreewise injective morphisms with degreewise projective cokernel.
\end{itemize}
We may define the full subcategory of non-positively graded $A_{\pal}$-modules $\PM^{\le0}(A_{\pal})\subseteq\PM(A_{\pal})$ simply replacing $\DGMod(A_{\alpha})$ by $\DGMod^{\le 0}(A_{\alpha})$. Notice that the same argument of Theorem~\ref{thm.modelpseudomodules} provides a model structure for $\PM^{\le0}(A_{\pal})$, where a morphism $\sF\to\sG$ in $\PM^{\le0}(A_{\pal})$ is a weak equivalence (respectively: cofibration, trivial fibration) if and only if it is a weak equivalence (respectively: cofibration, trivial fibration) in $\PM(A_{\pal})$. The same does not hold for fibrations. In particular, the natural inclusion functor $\PM^{\le 0}(A_{\pal})\to\PM(A_{\pal})$ is a left Quillen functor.
\end{remark}

There is a (homotopic) notion of quasi-coherent $A_{\pal}$-module.

\begin{definition}\label{def.qcohpseudo-module}
An $A_{\pal}$-module $\mathcal{F}$ over $X$ is called \textbf{quasi-coherent} if the morphism
\[ f_{\alpha\beta}\colon \sF_{\alpha}\otimes_{A_{\alpha}}A_{\beta} \to \sF_{\beta} \]
is a weak equivalence in $\DGMod(A_{\beta})$ for every $\alpha\leq\beta$ in $\sN$.
\end{definition}

We shall denote by $\QCoh(A_{\pal})\subseteq\PM(A_{\pal})$, and respectively by $\QCoh^{\ast}(A_{\pal})\subseteq\DGPM(A_{\pal})$, the full subcategories whose objects are quasi-coherent $A_{\pal}$-modules.
Every quasi-coherent sheaf over $X$ induces a quasi-coherent $A_{\pal}$-module in the obvious way.

\begin{remark}\label{rmk.qcoh-homotopyinvariance}
Quasi-coherent $A_{\pal}$-modules are closed under weak equivalences. Namely, given a weak equivalence $\varphi\colon\sF\to \sG$ in $\PM(A_{\pal})$ then $\sF$ is quasi-coherent if and only if $\sG$ is so. To prove the claim it is sufficient to consider the commutative diagram
\[ \xymatrix{	F_{\alpha}\otimes_{A_{\alpha}} A_{\beta} \ar@{->}[rr]^{f_{\alpha\beta}} \ar@{->}[d]^{\varphi_{\alpha}\otimes\Id} & & F_{\beta} \ar@{->}[d]^{\varphi_{\beta}} \\
G_{\alpha}\otimes_{A_{\alpha}}A_{\beta} \ar@{->}[rr]_{g_{\alpha\beta}} & & G_{\beta}	} \]
for every $\alpha\le \beta$ in $\sN$. The statement follows by the flatness of the map $A_{\alpha}\to A_{\beta}$ and by the \emph{2 out of 3} property.
\end{remark}

Observe that by Remark~\ref{rmk.qcoh-homotopyinvariance} the subcategory $\QCoh(A_{\pal})\subseteq\PM(A_{\pal})$ is closed under both factorizations given by Theorem~\ref{thm.modelpseudomodules}.

We now prove a preliminary result which will be useful later on.

\begin{lemma}\label{lemma.cofibrantprojective}
Let $\sQ\in\PM(A_{\pal})$ be a cofibrant $A_{\pal}$-module. Given a cospan $\sQ\xrightarrow{f} \sR\xleftarrow{\pi} \sP$ in $\DGPM(A_{\pal})$, if $\pi$ is degreewise surjective then there exists $h\in\Hom^{\ast}_{A_{\pal}}(\sQ,\sP)$ such that $\pi h=f$.
\begin{proof}
For simplicity we assume that $f\in\Hom^{0}_{A_{\pal}}(\sQ,\sR)$; the general case can be obtain by a shift. Fix $i\in\Z$; the map $\pi^i\colon \sR^i\to \sP^i$ induces the map of $A_{\pal}$-modules
\[ \xymatrix{	\hat{\sR} \colon \ar@{->}[d]^{\hat{\pi}} & \cdots \ar@{->}[r] & 0 \ar@{->}[r] \ar@{->}[d] & \sR^i \ar@{->}[r]^{\Id} \ar@{->}[d]^{\pi^i} & \sR^i \ar@{->}[r] \ar@{->}[d]^{\pi^i} & 0 \ar@{->}[r] \ar@{->}[d] & \cdots \\
\hat{\sP} \colon & \cdots \ar@{->}[r] & 0 \ar@{->}[r] & \sP^i \ar@{->}[r]^{\Id} & \sP^i \ar@{->}[r] & 0 \ar@{->}[r] & \cdots 	} \]
which is a trivial fibration.
Moreover, $f^{i}\colon \sQ^{i}\to \sR^{i}$ induces the map of $A_{\pal}$-modules
\[ \xymatrix{	\sQ \colon \ar@{->}[d]^{\hat{f}} & \cdots \ar@{->}[r] & \sQ^{i-2} \ar@{->}[r]^{d_{\sQ}^{i-2}} \ar@{->}[d] & \sQ^{i-1} \ar@{->}[r]^{d_{\sQ}^{i-1}} \ar@{->}[d]^{f^id_{\sQ}^{i-1}} & \sQ^i \ar@{->}[r] \ar@{->}[d]^{f^i} & \sQ^{i+1} \ar@{->}[r] \ar@{->}[d] & \cdots \\
\hat{\sP} \colon & \cdots \ar@{->}[r] & 0 \ar@{->}[r] & \sP^i \ar@{->}[r]^{\Id} & \sP^i \ar@{->}[r] & 0 \ar@{->}[r] & \cdots 	} \]
which can be lifted to $\hat{\sR}$ because $\sQ$ is cofibrant by assumption; i.e. there exists a map of $A_{\pal}$-modules $\hat{h}\colon\sQ\to \hat{\sR}$ such that $\hat{\pi}\hat{h}=\hat{f}$. Now define $h^i=\hat{h}^i\colon \sQ^i\to \sR^i$; reproducing the same argument for every $i\in\Z$ we obtain the required map $h\in\Hom_{A_{\pal}}^{0}(\sQ,\sP)$.
\end{proof}
\end{lemma}

Notice that if $X$ is an affine scheme then we can choose $\sN=\{\ast\}$. Therefore $A_{\pal}$-modules reduce to the category of DG-modules over $\Gamma(X,\Oh_X)$, and Lemma~\ref{lemma.cofibrantprojective} states that cofibrant DG-modules are degreewise projective. In the general case, the liftings $\{h_{\gamma}^i\colon Q_{\gamma}^i\to \sP^i_{\gamma}\}_{\gamma\in\sN}$ satisfy the commutativity relations induced by the nerve for any fixed $i\in\Z$.

\bigskip

\subsection{$A_{\pal}$-modules as sheaves over the nerve}\label{section.pseudosheaves}
Our next goal is to give a ``sheaf theoretic" description of $A_{\pal}$-modules. To this aim, we define a topology $\tau_{\sN}$ on the nerve $\sN$ as follows: $V\in\tau_{\sN}$ if and only if the condition
\[ \alpha\in V, \alpha\leq\beta \Rightarrow \beta\in V \]
is satisfied. This is called the Alexandroff topology, since $(\sN,\tau_{\sN})$ becomes an Alexandroff topological space, see~\cite{Alex}. For every fixed $\alpha\in\sN$ the set $V_{\alpha}=\{\gamma\in\sN\,\vert\,\alpha\leq\gamma\}\subseteq\sN$ is open, and the collection $\{V_{\alpha}\}_{\alpha\in\sN}\subseteq\tau_{\sN}$ is a basis for the topology. Then consider the category $\Sh_X(\sN)$ of sheaves of complexes over $\sN$; where moreover on every $V_{\alpha}$ it is given a structure of DG-module over $A_{\alpha}$ compatible with the restriction maps. Now, there is a pair of functors
\[ \sS\colon\PM(A_{\pal}) \to \Sh(\sN) \qquad \qquad \Gamma\colon \Sh(\sN) \to \PM(A_{\pal}) \]
defined by
\[ \sS(\sF)(V) =\lim_{\gamma\in V}\sF_{\gamma} \qquad \mbox{ and } \qquad \Gamma(\sG)_{\alpha} = \sG(V_{\alpha}) \]
for every $\sF\in\PM(A_{\pal})$, every $\sG\in\Sh(\sN)$, every $\alpha\in\sN$ and every $V\in\tau_{\sN}$. 
Notice that
\[ \sS(\sF)(V) = \left\{\{s_{\gamma}\}\in\prod_{\gamma\in V}\sF_{\gamma}\,\Bigl\vert\, f_{\gamma_1\gamma_2}(s_{\gamma_1}\otimes1)=s_{\gamma_2} \mbox{ for every }\gamma_1\leq\gamma_2 \right\} \]
and that $\sS(\sF)(V_{\alpha}) = \sF_{\alpha}$. for every $\alpha\in\sN$. In particular, $\Gamma\circ\sS=\Id_{\PM(A_{\pal})}$. Given $\sG\in\Sh_X(\sN)$ we have a natural map
\[ \sG(V) \xrightarrow{\cong} \lim_{\gamma\in V}\sG(V_{\gamma}) =  \sS(\Gamma(\sG))(V) \]
for every $V\in\tau_{\sN}$, which is an isomorphism because $\sG$ is a sheaf and $\bigcup_{\gamma\in V} V_{\gamma}=V$. Therefore the functors $\sS\colon\PM(A_{\pal}) \leftrightarrows \Sh(\sN)\colon\Gamma$ are equivalences of categories. A similar result can be found in~\cite[Proposition 6.6]{BBR}.

Recall that a sheaf $\sG$ of $\mathcal{O}_X$-modules is flasque if the restriction map $\sG(U)\to \sG(V)$ is surjective for every inclusion $V\to U$ between open subsets of $X$.

\begin{definition}\label{def.flasque}
An $A_{\pal}$-module $\sF\in\PM(A_{\pal})$ is called \textbf{flasque} if the associated sheaf $\sS(\sF)$ is so.
\end{definition}

\bigskip

\subsection{Inverse and direct image for $A_{\pal}$-modules: $j^{\ast}_{V}\dashv j_{V,\ast}$}
For any fixed open $V\in\tau_{\sN}$, denote by $j_V\colon V\hookrightarrow\sN$ the natural inclusion; the aim of this subsection is to introduce two functors $j_V^{\ast}$ and $j_{V,\ast}$, which we defined the ``inverse image" and ``direct image" functors because of the equivalence described in Subsection~\ref{section.pseudosheaves}.

First define $U_V=\bigcup_{\gamma\in V}U_{\gamma}\subseteq X$; recall that for every $\alpha\in\sN$ we denoted $V_{\alpha}=\{\gamma\in\sN\,\vert\, \gamma\geq\alpha\}$, so that in particular $U_{V_{\alpha}}=U_{\alpha}\subseteq X$.
Then the ``inverse image" and ``direct image" functors are defined by
\[ \begin{aligned}
j^{\ast}_V\colon \PM(A_{\pal}) &\to \PM(U_V)\\
\{ \sF_{\gamma} \}_{\gamma\in\sN} &\mapsto \{ \sF_{\gamma} \}_{\gamma\in V} 
\end{aligned} \qquad \mbox{ and } \qquad \begin{aligned}
j_{V,\ast}\colon \PM(U_V) &\to \PM(A_{\pal}) \\
\{ \sG_{\alpha} \}_{\alpha\in V} &\mapsto \left\{\lim\limits_{V\cap V_{\alpha}}\sG \right\}_{\alpha\in\sN}
\end{aligned} \]
respectively. More explicitly:
\[ (j_{V,\ast}\,\sG)_{\alpha} =\begin{cases}
\lim\limits_{\gamma\in V\cap V_{\alpha}}\sG_{\gamma} & \mbox{if } U_{\alpha}\cap U_{V}\neq\emptyset \\
0 & \mbox{otherwise}
\end{cases} \]
where the limit is taken in $\DGMod(A_{\alpha})$, and the $A_{\alpha}$-module structure is induced via $A_{\alpha}\to A_{\gamma}$ on each $\sG_{\gamma}$. 
Given $\alpha\leq\beta$ in $\sN$ such that $U_{\beta}\cap U_V\neq\emptyset$, the limit induces a natural map
\[ (j_{V,\ast}\,\sG)_{\alpha} = \lim_{\gamma\in V\cap V_{\alpha}}\sG_{\gamma} \longrightarrow \lim_{\gamma\in V\cap V_{\beta}}\sG_{\gamma} = (j_{V,\ast}\,\sG)_{\beta} \]
between DG-modules over $A_{\alpha}$. Since the $A_{\alpha}$ structure on $\lim\limits_{V\cap V_{\beta}}\sG_{\gamma}$ is given by $A_{\alpha}\to A_{\beta}$, by adjunction the above map corresponds to a morphism
\[ (j_{V,\ast}\,\sG)_{\alpha}\otimes_{A_{\alpha}}A_{\beta} \to (j_{V,\ast}\,\sG)_{\beta} \]
between DG-modules over $A_{\beta}$. Notice that in particular if $\alpha\in V$ then $(j_{V,\ast}\,\sG)_{\alpha}=\sG_{\alpha}$.

The following is a standard result, for which we provide an elementary proof in terms of the unit map since it will be useful later on.

\begin{lemma}\label{lemma.restriction}
For every open subset $j_V\colon V\hookrightarrow\sN$, there is an adjunction $j^{\ast}_{V}\dashv j_{V,\ast}$.
\begin{proof}
First notice that $j^{\ast}_{V}j_{V,\ast}$ is the identity on $\PM(U_{V})$. Hence to prove the adjunction it is sufficient to explicitly describe the unit $\eta\colon \Id_{\PM(A_{\pal})}\to j_{V,\ast}\, j^{\ast}_{V}$.
Fix $\sF\in\PM(A_{\pal})$ and $\gamma\in\sN$; the unit $\eta$ is defined on $\sF$ as
\[ \eta_{\sF} = \begin{cases}
\sF_{\alpha} \to \lim\limits_{\gamma\in V\cap V_{\alpha}}\sF_{\gamma}=(j_{V,\ast}\, j^{\ast}_{V}\sF)_{\alpha} & \mbox{if } U_{V}\cap U_{\alpha}\neq \emptyset \\
\sF_{\alpha} \to 0 & \mbox{otherwise}
\end{cases} \]
and the \emph{unit-counit equations} reduces to $\eta_{j_{V,\ast}\sG} = \Id_{j_{V,\ast}\sG}$ for every $\sG\in\PM(U_{V})$.
\end{proof}
\end{lemma}

\begin{remark}\label{rmk.restrictioncofibrant}
The adjoint pair of Lemma~\ref{lemma.restriction} is not necessarily a Quillen pair; in particular, the restriction $j_V^{\ast}\sF$ of a cofibrant $A_{\pal}$-module $\sF\in\PM(A_{\pal})$ may not be cofibrant. The crucial point is that the functor
\[ \lim_{V\cap V_{\alpha}}\colon \DGMod(A_{\alpha})^{V\cap V_{\alpha}}\to \DGMod(A_{\alpha}) \]
is right adjoint to the constant diagram, which does not preserve cofibrations in general. Nevertheless, if we choose $V=V_{\bar{\alpha}}=\{\gamma\in\sN\,\vert\,\bar{\alpha}\leq\gamma\}$ for some $\bar{\alpha}$ then the adjunction $j^{\ast}_{V_{\bar{\alpha}}}\dashv j_{V_{\bar{\alpha}},\ast}$ is in fact a Quillen pair. To prove the claim, notice that for every $\alpha\in\sN$ such that $U_{V_{\bar{\alpha}}}\cap U_{\alpha}\neq\emptyset$ we have $V_{\bar{\alpha}}\cap V_{\alpha} = V_{\bar{\alpha}\cup\alpha}$. Hence the constant functor
\[ \DGMod(A_{\alpha}) \to \DGMod(A_{\alpha})^{V_{\bar{\alpha}\cup\alpha}} \]
preserves cofibrations and trivial cofibrations; in fact for every $\beta\in\sN$ the set $\{\gamma\in\sN\,\vert\, \bar{\alpha}\cup\alpha\leq\gamma<\beta\}$ is connected. It follows that the functor $\lim_{V_{\bar{\alpha}\cup\alpha}}$ preserves fibrations and trivial fibrations, so that $j_{V_{\bar{\alpha}},\ast}\colon \PM(U_{V_{\bar{\alpha}}})\to \PM(A_{\pal})$ is a right Quillen functor as required. In particular, given a cofibrant $A_{\pal}$-module $\sF\in\PM(A_{\pal})$, its restriction $j_{V_{\bar{\alpha}}}^{\ast}\sF$ to $V_{\bar{\alpha}}$  is cofibrant in $\PM(U_{\bar{\alpha}})$.
\end{remark}

\begin{remark}
Notice that in the proof of Lemma~\ref{lemma.restriction} the differentials do not play any role. Therefore we have binatural isomorphisms
\[  \Hom_{A_{\pal},V}\left(j_{V}^{\ast}\sQ,\sG\right) \cong \Hom_{A_{\pal}}\left(\sQ,j_{V,\ast}\,\sG\right) \]
\[  \Hom^{\ast}_{A_{\pal},V}\left(j_{V}^{\ast}\sQ,\sG\right) \cong \Hom^{\ast}_{A_{\pal}}\left(\sQ,j_{V,\ast}\,\sG\right) \]
for every $\sQ\in\PM(A_{\pal})$ and every $\sG\in\PM(U_{\alpha})$. To avoid possible confusion we denoted morphisms in $\PM(U_{V})$ by $\Hom_{A_{\pal},V}(-,-)$, and $\ast$-morphisms in $\DGPM(U_{V})$ by $\Hom^{\ast}_{A_{\pal},V}(-,-)$.
\end{remark}

\begin{lemma}\label{lemma.restriction2}
Fix an open subset $j_V\colon V\hookrightarrow\sN$. Let $\sQ,\sG\in\PM(A_{\pal})$ and assume $\sQ$ to be cofibrant. Denote by $\eta_{\sG}\colon \sG\to j_{V,\ast}\,j_{V}^{\ast}\sG$ the unit map of the adjunction given by Lemma~\ref{lemma.restriction}. If $\eta_{\sG}$ is degreewise surjective, then the induced morphism
\[ \Hom_{A_{\pal}}^{\ast}(\sQ,\sG) \xrightarrow{ \eta_{\sG} } \Hom_{A_{\pal}}^{\ast}\left( \sQ, j_{V,\ast}\, j_{V}^{\ast}\sG \right)=\Hom_{A_{\pal},V}^{\ast}(j_{V}^{\ast}\sQ,j_{V}^{\ast}\sG) \]
is degreewise surjective.
\begin{proof}
We prove that the map $\Hom_{A_{\pal}}^{0}(\sQ,\sG) \xrightarrow{ \eta_{\sG} } \Hom_{A_{\pal},V}^{0}\left( \sQ, \Upsilon_{\ast}^{V}\Upsilon_{V}^{\ast}\sG \right)$ is surjective. The same argument works for other degrees. We need to show that every $\{\varphi_{\gamma}\}_{\gamma\in\sN}\in\Hom_{A_{\pal}}^{0}\left( \sQ, j_{V,\ast}\, j_{V}^{\ast}\sG \right)$ factors through the unit map $\eta_{\sG}$. Recall that since $\sQ$ is cofibrant then $\sQ^p$ is projective (in the sense of Lemma~\ref{lemma.cofibrantprojective}) for every $p\in\Z$, so that there exists the dotted morphism
\[ \xymatrix{	 & \sG^p \ar@{->}[d]^{\eta_{\sG}} \\
\sQ^p \ar@{->}[r] \ar@{.>}[ur] & (j_{V,\ast}\, j_{V}^{\ast}\sG)^p 	} \]
whence the statement.
\end{proof}
\end{lemma}

Lemma~\ref{lemma.restriction2} can be restated in terms of flasque $A_{\pal}$-modules, see Definition~\ref{def.flasque}. For every pair of $A_{\pal}$-modules $\sQ,\sG\in\PM(A_{\pal})$ it is defined an $A_{\pal}$-module $\Hom_{A_{\pal}}^{\ast}(\sQ,\sG)\in\PM(A_{\pal})$ as follows
\begin{enumerate}
\item $\Hom_{A_{\pal}}^{\ast}(\sQ,\sG)_{\alpha}=\Hom_{A_{\pal},V_{\alpha}}^{\ast}\left(j_{V_{\alpha}}^{\ast}\sQ,j_{V_{\alpha}}^{\ast}\sG\right) = \Hom^{\ast}_{A_{\pal}}\left(\sQ,j_{V_{\alpha},\ast}\, j_{V_{\alpha}}^{\ast}\sG\right)$ for every $\alpha\in\sN$,
\item \[ \begin{aligned}
\Hom_{A_{\pal}}^{\ast}(\sQ,\sG)_{\alpha} \otimes_{A_{\alpha}}A_{\beta} &\to \Hom_{A_{\pal}}^{\ast}(\sQ,\sG)_{\beta} \\
\{\varphi_{\gamma}\}_{\gamma\geq\alpha}\otimes x &\mapsto \{x\cdot\varphi_{\gamma}\}_{\gamma\geq\beta}
\end{aligned} \]
for every $\alpha\leq\beta$ in $\sN$.
\end{enumerate}

\begin{proposition}
Let $\sQ,\sG\in\PM(A_{\pal})$; assume $\sQ$ to be cofibrant and $\sG$ to be flasque. Then the $A_{\pal}$-module $\Hom_{A_{\pal}}^{\ast}(\sQ,\sG)\in\PM(A_{\pal})$ is flasque.
\begin{proof}
If $\sG$ is flasque then for every open subset $j_V\colon V\hookrightarrow \sN$ the unit map $\eta_{\sG}\colon \sG\to j_{V,\ast}\, j_{V}^{\ast}\sG$ given by Lemma~\ref{lemma.restriction} is surjective. The statement follows by Lemma~\ref{lemma.restriction2}.
\end{proof}
\end{proposition}

\bigskip

\section{Extended lower-shriek functor}\label{section.lower-shriek}

This section is devoted to the well posedness of a certain functor that we shall call the \emph{extended lower-shriek}.
Again, $X$ is a finite-dimensional separated Noetherian scheme over a field $\K$ of characteristic $0$. Moreover, a fixed affine open cover $\{U_h\}_{h\in H}$ is chosen, and we denote its nerve by $\sN$.

\begin{definition}\label{def.L}
Define the poset $\LI$ as
\begin{enumerate}
\item $\LI= \{(\beta,\gamma)\in \sN\times \sN \,\vert\, \beta\leq\gamma\}$,
\item  $(\beta,\gamma)\le (\delta,\eta)$ if and only if $\beta\leq\delta$ and $\eta\leq\gamma$ in $\sN$.
\end{enumerate}
\end{definition}

In particular, condition $(2)$ of Definition~\ref{def.L} implies that for every $\beta\leq\delta\leq\eta\leq\gamma$ the diagram
\[ \xymatrix{	(\beta,\gamma) \ar@{->}[r] \ar@{->}[d] \ar@{->}[dr] & (\delta,\gamma) \ar@{->}[d] \\
(\beta,\eta) \ar@{->}[r] & (\delta,\eta)	} \]
commutes in $\LI$. We shall call a morphism $(\beta,\gamma)\to (\delta,\gamma)$ an \emph{horizontal} morphism, and similarly we call morphisms of the form $(\beta,\gamma)\to (\beta,\eta)$ \emph{vertical} morphisms.

\begin{remark}
More generally, for every small category $C$ we can consider the category 
$\mathbf{L}_C$ whose objects are maps in $C$ and whose morphisms are commutative diagrams: 
\[ \xymatrix{	\beta\ar[r] \ar[d]& \gamma &\ar@{}[d]^{\iff}&(\beta\to\gamma)\ar[d]^{\in\Mor_{\mathbf{L}_C}}\\
\delta\ar[r] & \eta\ar[u]&&(\delta\to\eta)} \]
If $C$ is a direct Reedy category, then  $\mathbf{L}_C$ is an inverse Reedy category with degree function 
\[ \deg(\beta\to\gamma)=\deg(\gamma)-\deg(\beta)\ge 0\,.\]
\end{remark}

For every $\alpha\le \beta$ in $\sN$ denote by 
\[ i_{\beta}\colon U_{\beta}\xrightarrow{i_{\beta}^{\alpha}}U_{\alpha}\xrightarrow{i_{\alpha}}X\]
the natural inclusions.

Since the scheme is separated, then $U_{\alpha}$ is affine for every $\alpha\in\sN$. Hence the datum of an $A_{\pal}$-module $\sF\in\PM(A_{\pal})$ is equivalent to $\sF_\alpha\in \DGMod(\mathcal{O}_{U_\alpha})$ for every $\alpha\in\sN$ and morphisms
\[ f_{\alpha\beta} \colon {i_{\beta}^{\alpha}}^{\ast}\sF_{\alpha}=\sF_{\alpha}\vert_{U_\beta}\to \sF_{\beta},\qquad \alpha\le\beta\,.\]
Now, we fix the $A_{\pal}$-module $\sF$ and define the following functors
\[ \begin{aligned}
\mathcal{F}_{\ast}\colon \LI &\to \DGMod(\mathcal{O}_X) \qquad \qquad \sF_{!}\colon \LI \to \DGMod(\mathcal{O}_X) \\
(\beta,\gamma) &\mapsto {i_\gamma}_{\ast}{i_{\gamma}^{\beta}}^{\ast}\sF_{\beta}={i_\gamma}_{\ast}\sF_{\beta}\vert_{U_{\gamma}} \qquad 
(\beta,\gamma)\mapsto {i_\gamma}_!{i_{\gamma}^{\beta}}^{\ast}\sF_{\beta}={i_\gamma}_!(\sF_{\beta}\vert_{U_{\gamma}}) \; .
\end{aligned} \] 

If $(\beta,\gamma)\to (\delta,\eta)$ then $U_\gamma\subset U_{\eta}\subset U_{\delta}\subset U_{\beta}$, so that it is given the map
$f_{\beta\delta} \colon \sF_{\beta}\vert_{U_\delta}\to \sF_{\delta}$ which in turn induces the morphism $\mathcal{F}_{!}(\beta,\gamma)\to \mathcal{F}_{!}(\delta,\eta)$ defined by the composition
\[ {i_{\gamma}}_!(\sF_{\beta}\vert_{U_\gamma})\xrightarrow{{i_{\gamma}}_! (f_{\beta\delta}\vert_{U_{\gamma}})}{i_{\gamma}}_!\sF_{\delta}\vert_{U_\gamma}
\to  {i_{\eta}}_!\sF_{\delta}\vert_{U_\eta} \;.\]
Similarly, the morphisms  $\mathcal{F}_{{\ast}}(\beta,\gamma)\to \mathcal{F}_{{\ast}}(\delta,\eta)$ is given by the composition
\[ {i_{\gamma}}_{\ast}(\sF_{\beta}\vert_{U_\gamma})\xrightarrow{{i_{\gamma}}_{\ast} (f_{\beta\delta}\vert_{U_{\gamma}})}{i_{\gamma}}_{\ast}\sF_{\delta}\vert_{U_\gamma}
\to  {i_{\eta}}_{\ast}{\sF_{\delta}}_{|U_\eta} \;. \]
 
\begin{definition}\label{def.lowershriek}
In the above notation, the \textbf{extended lower-shriek} functor $\Upsilon_!$ is defined as
\[ \begin{aligned}
\Upsilon_! \colon \PM(A_{\pal}) &\to \DGMod(\mathcal{O}_X) \\
\mathcal{F} &\mapsto \colim\limits_{\LI} \sF_! \; .
\end{aligned} \]
\end{definition}

Our goal is now to investigate the relation between the extended lower-shriek and the following functor
\[ \begin{aligned}
\Upsilon^{\ast} \colon \DGMod(\mathcal{O}_X) &\to \PM(A_{\pal}) \\
\mathcal{F} &\mapsto \left\{ \mathcal{F}(U_{\alpha})\right\}_{\alpha\in\sN} \; .
\end{aligned} \]

\begin{proposition}\label{prop.shriekadjoint}
The functors $\Upsilon_! \colon \PM(A_{\pal}) \leftrightarrows \DGMod(\mathcal{O}_X)\colon\Upsilon^{\ast}$ form an adjoint pair.
\begin{proof}
We need to show that there exists a bi-natural bijection of sets
\[ \Hom_{\DGMod(\mathcal{O}_X)}(\Upsilon_!\mathcal{F},\mathcal{G}) \cong \Hom_{A_{\pal}}(\mathcal{F},\Upsilon^{\ast}\mathcal{G}) \]
for every $\mathcal{F}\in\PM(A_{\pal})$ and every $\mathcal{G}\in\DGMod(\mathcal{O}_X)$.
By the universal property of the colimit, the data of a morphism $\varphi\in\Hom_{\DGMod(\mathcal{O}_X)}(\Upsilon_!\mathcal{F},\mathcal{G})$ is equivalent to the following chain of one-to-one correspondences
\[ \begin{aligned}
\varphi & \longleftrightarrow \left\{i_{\gamma!}\left(\sF_{\beta}\vert_{U_{\gamma}}\right) \to \sG\right\}_{(\beta,\gamma)\in\LI} \longleftrightarrow \left\{\left(\sF_{\beta}\vert_{U_{\gamma}}\right) \to \sG\vert_{U_{\gamma}}\right\}_{(\beta,\gamma)\in\LI} \stackrel{(\ast)}{\longleftrightarrow} \\
& \stackrel{(\ast)}{\longleftrightarrow} \left\{ \sF_{\beta}(U_{\beta})\otimes_{A_{\beta}}A_{\gamma} \to \sG(U_{\gamma})\right\}_{(\beta,\gamma)\in\LI} \stackrel{(\ast\ast)}{\longleftrightarrow} \{\sF_{\gamma}(U_{\gamma})\to \sG(U_{\gamma})\}_{\gamma\in\sN}\in \Hom_{A_{\pal}}(\sF,\Upsilon^{\ast}\sG)
\end{aligned} \]
where:
\begin{itemize}
\item $(\ast)$ is a bijection since the morphisms of sheaves are all determined by localizations of the module $\sF_{\beta}\otimes_{A_{\beta}}A_{\gamma}$,
\item $(\ast\ast)$ is a bijection since for every $(\beta,\gamma)\in\LI$ we have a commutative diagram
\[ \xymatrix{	\sF_{\beta}(U_{\beta})\otimes_{A_{\beta}}A_{\gamma}\ar@{->}[r]^-{f_{\beta\gamma}} \ar@{->}[dr] & \sF_{\gamma}(U_{\gamma})\ar@{->}[d]\\
 & \sG(U_{\gamma})	} \]
where the morphisms $f_{\beta\gamma}$ are given by the $A_{\pal}$-module $\sF$.
\end{itemize}
\end{proof}
\end{proposition}

Recall that an object $\sF\in\DGMod(\mathcal{O}_X)$ is called a \textbf{flasque complex} if it is degreewise flasque, see~\cite{Hov01}.

\begin{theorem}\cite[Theorem 5.2]{Hov01}\label{thm.flatmodelstructure}
Let $X$ be a separated finite-dimensional Noetherian scheme. Then the category $\DGMod(\mathcal{O}_X)$ is endowed with the \textbf{flat model structure}, where the weak equivalences are the quasi-isomorphisms, and fibrations are epimorphisms with flasque kernel.
\end{theorem}

\begin{remark}\cite[Exercise II.1.6]{Har}\label{rmk.flatfibrations}
Let $\varphi\colon \mathcal{F}\to\mathcal{G}$ be an epimorphism of sheaves of $\mathcal{O}_X$-modules with flasque kernel over a separated Noetherian scheme $X$. Then $\varphi_V\colon\mathcal{F}(V)\to\mathcal{G}(V)$ is surjective for every open subset $V\subseteq X$.
\end{remark}

Our next result is a refined version of Proposition~\ref{prop.shriekadjoint}.

\begin{theorem}\label{thm.lowershriekQuillen}
The adjoint functors
\[ \Upsilon_! \colon \PM(A_{\pal}) \rightleftarrows \DGMod(\mathcal{O}_X)\colon\Upsilon^{\ast} \]
form a Quillen pair with respect to the model structure of Theorem~\ref{thm.modelpseudomodules} on $\PM(A_{\pal})$, and the flat model structure on $\DGMod(\mathcal{O}_X)$.
\begin{proof}
The adjointness follows from Proposition~\ref{prop.shriekadjoint}, and the right adjoint $\Upsilon^{\ast}$ preserves fibrations by Remark~\ref{rmk.flatfibrations}. In order to prove that the functor $\Upsilon^{\ast}$ preserves trivial fibrations it is sufficient to observe that the complex of sections $\Gamma(V,\ker(f))$ is acyclic for every open $V\subseteq X$ and for any epimorphism with flasque kernel $f\colon \sF\to\sG$ in $\DGMod(\Oh_X)$; this immediately follows from~\cite[Lemma 4.1]{Hov01}.
\end{proof}
\end{theorem}

Notice that the proof of Theorem~\ref{thm.lowershriekQuillen} relies on~\cite[Lemma 4.1]{Hov01}, which applies because we assumed the scheme $X$ to be Noetherian and finite-dimensional.

As an immediate consequence of Theorem~\ref{thm.lowershriekQuillen}, we obtain the existence of the total derived functors
\[ \mathbb{L}\Upsilon_! \colon \Ho(\PM(A_{\pal})) \leftrightarrows \Ho(\DGMod(\mathcal{O}_X))\colon \R\Upsilon^{\ast} \; . \]

\bigskip

\section{From $A_{\pal}$-modules to derived categories}\label{section.derivedlowershriek}

The first goal of this section is to show that the total left derived functor of the extended lower-shriek introduced in the Section~\ref{section.lower-shriek} maps (classes of) quasi-coherent $A_{\pal}$-modules in (classes of) complexes of quasi-coherent sheaves, see Theorem~\ref{thm.leftderivedlowershriek}. Hence there will be induced functors
\[ \overline{\mathbb{L}\Upsilon_!} \colon \Ho(\QCoh(A_{\pal})) \leftrightarrows \D(\QCoh(X)) \colon \overline{\R\Upsilon}^{\ast} \; . \]
Our main result shows that the above functors are in fact equivalences of triangulated categories, see Theorem~\ref{thm.equivalence}. To this aim, we shall first prove that
\[ \overline{\mathbb{L}\Upsilon_!}[\sF] = [\Upsilon_!\sF] \qquad \text{ for every } [\sF]\in\Ho(\QCoh(A_{\pal})) \]
\[ \overline{\R\Upsilon}^{\ast}[\sG] = [\Upsilon^{\ast}\sG] \qquad \text{ for every } [\sG]\in\D(\QCoh(X)) \; . \]

As usual, $X$ is a fixed separated finite-dimensional Noetherian scheme over $\K$; moreover $\sN$ denotes the nerve of a fixed affine open covering $\{U_h\}_{h\in H}$.
Recall that by Definition~\ref{def.qcohpseudo-module}, an $A_{\pal}$-module $\sF\in\PM(A_{\pal})$ is called \emph{quasi-coherent} if the morphism
\[ f_{\alpha\beta}\colon \sF_{\alpha}\otimes_{A_{\alpha}}A_{\beta} \to \sF_{\beta} \]
is a weak equivalence (i.e. a quasi-isomorphism) in $\DGMod(A_{\beta})$ for every $\alpha\leq\beta$ in $\sN$.

We need an easy preliminary result.

\begin{lemma}\label{lem.colimit-cohomology-commute}
Let $N$ be a small direct category and let $R$ be a ring. Consider the category $\DGMod(R)$ of complexes of $R$-modules. Given a functor $F\colon N\to\DGMod(R)$ there exists a natural isomorphism of $R$-modules $H^j\left(\colim_{\beta\in N}F_{\beta}\right) \cong \colim_{\beta\in N}(H^j(F_{\beta}))$ for every $j\in\Z$.
\begin{proof}
Consider the exact sequence $0\to Z^jF_{\beta}\to F_{\beta}^j\stackrel{d_{\beta}^j}{\longrightarrow} Z^{j+1}F_{\beta}\to H^{j+1}F_{\beta}\to 0$, for every $\beta\in N$ and every $j\in\Z$.
Now observe that the functor $\colim_{ N}$ is exact, being direct on a category of modules.
In particular,
\[\colim_{\beta\in N}Z^jF_{\beta} \cong\ker\left\{\colim_{\beta\in N}d_{\beta}^j\right\} = Z^j\left(\colim_{\beta\in N}F_{\beta}\right),\]
and the thesis easily follows.
\end{proof}
\end{lemma}

\begin{proposition}\label{prop.qiso-qcoherent}
Let $\sF\in\QCoh(A_{\pal})$ be a quasi-coherent $A_{\pal}$-module.
Then for every $\alpha\in \sN$ there exists a quasi-isomorphism $\widetilde{\sF_{\alpha}} \to (\Upsilon_!\mathcal{F})\vert_{U_{\alpha}}$ in $\DGMod(\Oh_{U_{\alpha}})$.
\begin{proof}
We show that the natural morphism
\[ \varphi\colon \widetilde{\sF}_{\alpha} \to \left(\colim_{(\beta,\gamma)\in\LI}i_{\gamma!}(\widetilde{\sF}_{\beta}\vert_{U_{\gamma}})\right)\Bigl\vert_{U_{\alpha}} = (\Upsilon_!\sF)\vert_{U_{\alpha}}\]
is a quasi-isomorphism by showing that the induced morphism $\varphi_{x}$ is so at each stalk, ${x}\in U_{\alpha}$. Consider the following chain of equalities
\[ \left((\Upsilon_!\mathcal{F})\vert_{U_{\alpha}}\right)_{x} = \colim_{(\beta,\gamma)\in\LI}\left(i_{\gamma!}(\widetilde{\sF}_{\beta}\vert_{U_{\gamma}})\right)_{x} = \colim_{\{(\beta,\gamma)\in\LI\,\vert\, x\in U_{\gamma}\}} \left(\widetilde{\sF}_{\beta}\vert_{U_{\gamma}}\right)_{x} = \colim_{\beta\in\mathcal{N}} \left(\widetilde{F}_{\beta}\right)_{x} \]
where the last equality holds since for every $\beta\leq\gamma_1\leq\gamma_2$ the vertical morphism induced on the stalk $\left(\widetilde{\sF}_{\beta}\vert_{U_{\gamma_1}}\right)_{x} \to \left(\widetilde{\sF}_{\beta}\vert_{U_{\gamma_2}}\right)_{x}$ is an isomorphism, being $x\in U_{\gamma_2}\subseteq U_{\gamma_1}$. Now take $j\in\mathbb{Z}$ and notice that $\mathcal{N}$ is connected, whenever $\beta_1\leq\beta_2$ the natural morphism $H^j(\widetilde{F}_{\beta_1})_{x}\to H^j(\widetilde{F}_{\beta_2})_{x}$ is an isomorphism by hypothesis; hence
\[H^j(\varphi_{x})\colon H^j(\widetilde{F}_{\alpha})_{x}\stackrel{\cong}{\longrightarrow} \colim_{\beta\in\mathcal{N}}H^j(\widetilde{F}_{\beta})_{x} \cong [\mbox{Lemma } \ref{lem.colimit-cohomology-commute}] \cong H^j\left(\colim_{\beta\in\mathcal{N}}(\widetilde{F}_{\beta})\right)_{x}\]
and the statement follows.
\end{proof}
\end{proposition}

Notice that there are inclusion functors
\[ \Ho(\QCoh(A_{\pal})) \to \Ho(\PM(A_{\pal})) \mbox{ and } \D(\QCoh(X))\to \Ho(\DGMod(\mathcal{O}_X)) \, . \]
Our goal is now to show that the total left derived functor $\mathbb{L}\Upsilon_! \colon \Ho(\PM(A_{\pal})) \to \Ho(\DGMod(\mathcal{O}_X))$ maps $\Ho(\QCoh(A_{\pal}))$ to $\D(\QCoh(X))$.

\begin{remark}\label{rmk.neeman}
Let $\D_{qc}\left(\mathcal{O}_X\right)$ be the derived category of cochain complexes of arbitrary $\mathcal{O}_X$-modules over $X$, with quasi-coherent cohomology. Then the natural functor $\D(\QCoh(X))\to \D_{qc}\left(\mathcal{O}_X\right)$ is an equivalence of categories, see~\cite{BN}.
\end{remark}

\begin{theorem}\label{thm.leftderivedlowershriek}
The functor $\mathbb{L}\Upsilon_! \colon \Ho(\PM(A_{\pal})) \to \Ho(\DGMod(\mathcal{O}_X))$ maps (classes of) quasi-coherent $A_{\pal}$-modules to (classes of) complexes of quasi-coherent sheaves.
\begin{proof}
The statement immediately follows by Proposition~\ref{prop.qiso-qcoherent} and Remark~\ref{rmk.neeman}.
\end{proof}
\end{theorem}

The functor $\Upsilon^{\ast}$ obviously maps quasi-coherent sheaves to quasi-coherent $A_{\pal}$-modules. Therefore by Theorem~\ref{thm.leftderivedlowershriek} the restricted functors
\[ \overline{\mathbb{L}\Upsilon_!} \colon \Ho(\QCoh(A_{\pal})) \leftrightarrows \D(\QCoh(X))\colon \overline{\R\Upsilon^{\ast}} \]
are well-defined.

\bigskip

\subsection{The equivalence $\Ho(\QCoh(A_{\pal}))\simeq\D(\QCoh(X))$}\label{section.equivalence}

The aim of this subsection is to show that the adjoint pair
\[ \overline{\mathbb{L}\Upsilon}_! \colon \Ho(\QCoh(A_{\pal})) \leftrightarrows \D(\QCoh(X))\colon \overline{\mathbb{R}\Upsilon}^{\ast} \]
introduced in the section above is in fact an equivalence of triangulated categories.

Explicit models for the (unique) DG-enhancement of $\D(\QCoh(X))$ already exist, e.g. the category of complexes of injectives. For a survey concerning this topic we refer to~\cite{CS} and~\cite{LO}. As we shall see, cofibrant $A_{\pal}$-modules produce another explicit DG-enhancement for $\D(\QCoh(X))$, see Corollary~\ref{corollary.enhancement}.

\begin{remark}\label{rmk.creator}
The functor $\Upsilon^{\ast}\colon \DGMod(\Oh_X)\to \PM(A_{\pal})$ maps quasi-isomorphisms between (complexes of) quasi-coherent sheaves to weak equivalences between quasi-coherent $A_{\pal}$-modules. This easily follows recalling that cohomology commutes with direct colimits (hence with stalks), see Lemma~\ref{lem.colimit-cohomology-commute}.
In particular, $\overline{\mathbb{R}\Upsilon}^{\ast}[\sF]=[\Upsilon^{\ast}(\sF)]$ for every $[\sF]\in \D(\QCoh(X))$.
\end{remark}

\begin{lemma}\label{lemma.counit}
Let $\varphi\colon\sF\to \sG$ be a morphism in $\QCoh(A_{\pal})$. Then $\varphi$ is a weak equivalence if and only if $\Upsilon_!(\varphi)$ is a weak equivalence in $\DGMod(\Oh_X)$.
\begin{proof}
For any $\alpha\in\sN$ consider the commutative diagram
\[ \xymatrix{	\widetilde{\sF_{\alpha}} \ar@{->}[rr] \ar@{->}[d] & & \Upsilon_!(\sF)\vert_{U_{\alpha}} \ar@{->}[d] \\
\widetilde{\sG_{\alpha}} \ar@{->}[rr] & & \Upsilon_!(\sG)\vert_{U_{\alpha}}		} \]
where the horizontal arrows are quasi-isomorphisms in $\DGMod(\Oh_{U_{\alpha}})$ by Proposition~\ref{prop.qiso-qcoherent}. Observe that $\sF_{\alpha}\to\sG_{\alpha}$ is a quasi-isomorphism in $\DGMod(A_{\alpha})$ if and only if $\widetilde{\sF_{\alpha}}\to\widetilde{\sG_{\alpha}}$ is so on each stalk in $U_{\alpha}$. Then the statement follows by the \emph{2 out of 3} property.
\end{proof}
\end{lemma}

Notice that Lemma~\ref{lemma.counit} implies that $\overline{\mathbb{L}\Upsilon}_![\sG] = [\Upsilon_!\sG]$ for every $[\sG]\in\Ho(\QCoh(A_{\pal}))$. Hence it is convenient to simply denote by
\[ \Upsilon_! \colon \Ho(\QCoh(A_{\pal})) \leftrightarrows \D(\QCoh(X))\colon \Upsilon^{\ast} \]
the functors $\overline{\mathbb{L}\Upsilon}_!$ and $\overline{\mathbb{R}\Upsilon}^{\ast}$.

\begin{theorem}\label{thm.equivalence}
The functors $\Upsilon_! \colon \Ho(\QCoh(A_{\pal})) \leftrightarrows \D(\QCoh(X))\colon \Upsilon^{\ast}$ are equivalences of triangulated categories.
\begin{proof}
In order to avoid possible confusion, throughout all the proof we shall keep the notation $\overline{\mathbb{L}\Upsilon}_!$ and $\overline{\mathbb{R}\Upsilon}^{\ast}$ to denote the functors in the statement.

First recall that the triangulated structure is preserved because the functors come from a Quillen adjunction.
Hence we only need to prove that the natural morphisms
\[ \overline{\mathbb{L}\Upsilon}_!\circ \overline{\mathbb{R}\Upsilon}^{\ast}[\sF] \to [\sF] \qquad \mbox{ and } \qquad [\sG]\to \overline{\mathbb{R}\Upsilon}^{\ast}\circ\overline{\mathbb{L}\Upsilon}_![\sG] \]
are isomorphisms for every $[\sF]\in \D(\QCoh(X))$ and every $[\sG]\in\Ho(\PM(A_{\pal}))$.

\begin{enumerate}
\item First observe that $\overline{\mathbb{L}\Upsilon}_!\circ \overline{\mathbb{R}\Upsilon}^{\ast}[\sF]  = [\Upsilon_!\Upsilon^{\ast}(\sF)]$ by Remark~\ref{rmk.creator} and Lemma~\ref{lemma.counit}. Moreover, since
\[ \left( \Upsilon_!\Upsilon^{\ast}(\sF) \right)_{x} = \colim\limits_{(\beta,\gamma)\in\LI}(i_{\gamma!}(\sF\vert_{U_\gamma}))_{x} = \colim\limits_{\{(\beta,\gamma)\in\LI \,\vert\, x\in U_{\gamma}\}} (i_{\gamma!}(\sF\vert_{U_{\gamma}}))_{x} = \colim\limits_{\beta\in I} (\sF\vert_{U_{\beta}})_{x} = \sF_{x}	\]
for every $x\in X$, then the natural map $\Upsilon_!\Upsilon^{\ast}(\sF)\to\sF$ is an isomorphism.
\item The second natural isomorphism follows by Lemma~\ref{lemma.counit} and Proposition~\ref{prop.qiso-qcoherent}.
\end{enumerate}
\end{proof}
\end{theorem}

Theorem~\ref{thm.equivalence} partially appears in~\cite[Proposition 2.28]{BF}, where it is proven that $\Upsilon^{\ast}$ is an equivalence on its image.

Define the DG-category $\QCoh^{\ast}(A_{\pal})^{c}$ whose objects are cofibrant quasi-coherent $A_{\pal}$-modules, and whose morphisms are $\ast$-morphisms, see Definition~\ref{definition.DGpseudomap}.
Notice that
\[ Z^0\left(\QCoh^{\ast}(A_{\pal})^{c}\right) = \QCoh(A_{\pal})^{c} \; . \]
Moreover, every weak equivalence $\sF\to\sG$ in $\PM(A_{\pal})$ between cofibrant $A_{\pal}$-modules is in fact an isomorphism up to homotopy; i.e. $H^0\left(\QCoh^{\ast}(A_{\pal})^{c}\right) \simeq \Ho\left(\QCoh(A_{\pal})^{c}\right)$.

\begin{corollary}\label{corollary.enhancement}
The DG-category $\QCoh^{\ast}(A_{\pal})^{c}$ is a DG-enhancement for the unbounded derived category $\D(\QCoh(X))$.
\begin{proof}
There are equivalences of triangulated categories
\[ H^0\left(\QCoh^{\ast}(A_{\pal})^{c}\right) \simeq \Ho\left(\QCoh(A_{\pal})^{c}\right) \simeq \Ho\left(\QCoh(A_{\pal})\right)\simeq\D(\QCoh(X)) \;,\]
where the last one follows by Theorem~\ref{thm.equivalence}.
\end{proof}
\end{corollary}

\bigskip

\section{Derived endomorphisms of quasi-coherent sheaves}\label{section.REnd}

Throughout this section we shall consider a fixed finite-dimensional Noetherian separated scheme $X$ over a field $\K$, together with a quasi-coherent sheaf $\sF$ on it. Also, we fix an open affine covering $\{U_h\}_{h\in H}$, denoting by $\sN$ its nerve.

The first main goal of this section is to give different constructions of the derived endomorphisms $\R\End(\sF)$. The interest in this object arises in several areas of Algebraic Geometry; for instance it carries a DG-Lie structure controlling infinitesimal deformations of $\sF$ as we shall see in Section~\ref{section.deformations}.

Recall that $\R\End(\sF)$ is represented (up to quasi-isomorphisms) by the complex $\Hom^{\ast}_{\Oh_X}(\sF,\sI)$, for any injective resolution $\sF\to\sI$. Notice that $\Hom^{\ast}_{\Oh_X}(\sF,\sI)=\Hom_{\Oh_X}(\sF,\sI)$, up to a sign on the differential.

\bigskip

\subsection{$\R\End(\sF)$ via $A_{\pal}$-modules}

The aim of this subsection is to prove that given a cofibrant replacement $\varepsilon\colon \sQ\to \Upsilon^{\ast}\sF$ in $\PM(A_{\pal})$, then the derived endomorphisms of $\sF$ are represented by $\End_{A_{\pal}}^{\ast}(\sQ)$.

For notational convenience we shall also denote by $\varepsilon$ the induced map $\Upsilon_!\sQ\to\Upsilon_!\Upsilon^{\ast}\sF = \sF$.

\begin{proposition}\label{prop.REnd}
Let $\sF$ be a quasi-coherent sheaf on $X$, and consider a cofibrant replacement $\varepsilon\colon\sQ\to \Upsilon^{\ast}\sF$ in $\PM(A_{\pal})$. Then the induced map
\[ \Hom^{\ast}_{\Oh_X}(\Upsilon_!\sQ,\sJ) \xleftarrow{-\circ\varepsilon} \Hom^{\ast}_{\Oh_X}(\sF,\sJ)\]
is a quasi-isomorphism for any bounded below complex of injectives $\sJ$.
\begin{proof}
Since $\sJ$ is degreewise injective we have a short exact sequence
\[ 0\to \Hom^{\ast}_{\Oh_X}(\sF,\sJ)\to \Hom_{\Oh_X}^{\ast}(\Upsilon_!\sQ,\sJ) \to \Hom^{\ast}_{\Oh_X}(\sH,\sJ)\to 0 \]
where $\sH=\ker(\varepsilon)$ is acyclic. By standard arguments it is easy to show that any map from an acyclic complex to a bounded below complex of injectives is homotopic to zero, see e.g.~\cite[III.5.24]{GelfandManin}. Hence the complex $\Hom^{\ast}_{\Oh_X}(\sH,\sJ)$ is acyclic and the statement follows.
\end{proof}
\end{proposition}

\begin{proposition}\label{proposition.REnd2}
Let $\sF$ be a quasi-coherent sheaf on $X$, let $\varphi\colon\sF\to \sI$ be an injective resolution, and consider a cofibrant replacement $\varepsilon\colon\sQ\to \Upsilon^{\ast}\sF$ in $\PM(A_{\pal})$. Then the maps
\[ \Hom^{\ast}_{A_{\pal}}(\sQ,\sQ) \xrightarrow{\varepsilon\circ-}\Hom^{\ast}_{A_{\pal}}(\sQ,\Upsilon^{\ast}\sF) \xrightarrow{\varphi\circ-} \Hom^{\ast}_{A_{\pal}}(\sQ,\Upsilon^{\ast}\sI) \]
are quasi-isomorphisms.
\begin{proof}
We shall prove that the functor $\Hom_{A_{\pal}}^{\ast}(\sQ,-)\colon \PM(A_{\pal})\to \DGMod(\Z)$ maps weak equivalences to quasi-isomorphisms, being $\sQ$ cofibrant.
Since every object in $\PM(A_{\pal})$ is fibrant, by Ken Brown's Lemma it is sufficient to show that $\Hom_{A_{\pal}}^{\ast}(\sQ,-)$ maps trivial fibrations to quasi-isomorphisms. To this aim, take a trivial fibration $f\colon \sG\to \sH$ in $\PM(A_{\pal})$. Then we have a short exact sequence
\[ 0 \to \Hom^{\ast}_{A_{\pal}}(\sQ,\ker(f)) \to \Hom^{\ast}_{A_{\pal}}(\sQ,\sG) \xrightarrow{f\circ-} \Hom_{A_{\pal}}^{\ast}(\sQ,\sH) \to 0 \; ;\]
where the surjectivity comes from Lemma~\ref{lemma.cofibrantprojective}.

To conclude we need to show that $\Hom^{\ast}_{A_{\pal}}(\sQ,\ker(f))$ is acyclic. Notice that every cocycle $[h]\in Z^n\left(\Hom^{\ast}_{A_{\pal}}(\sQ,\ker(f))\right)$ is given by a map $h\colon \sQ\to \ker(f)[n]$ of $A_{\pal}$-modules. Now, factor the weak equivalence $0\to \ker(f)$ as
\[ 0 \xrightarrow{\iota} \cocone\left(\Id_{\ker(f)[n]}\right) \xrightarrow{\pi} \ker(f)[n]  \]
and observe that $\iota$ is a weak equivalence and $\pi$ is a trivial fibration. Hence the square of solid arrows
\[ \xymatrix{	0 \ar@{->}[r]^-{\iota} \ar@{->}[d] & \cocone\left(\Id_{\ker(f)[n]}\right) \ar@{->}[d]^{\pi} \\
\sQ \ar@{->}[r]_-{h} \ar@{.>}[ur]_{\bar{h}} & \ker(f)[n]	} \]
admits the dotted lifting $\bar{h}\colon \sQ\to \cocone\left(\Id_{\ker(f)[n]}\right)$, which in turn implies that $h$ is homotopic to zero, i.e. $[h]=[0]\in H^n\left(\Hom^{\ast}_{A_{\pal}}(\sQ,\ker(f))\right)$.
\end{proof}
\end{proposition}

\begin{remark}\label{rmk.REnd2}
The same argument given in the proof of Proposition~\ref{proposition.REnd2} leads to quasi-isomorphisms
\[ \Hom^{\ast}_{A_{\pal}}(\sQ,\sQ)_{\alpha} \xrightarrow{\varepsilon\circ-}\Hom^{\ast}_{A_{\pal}}(\sQ,\Upsilon^{\ast}\sF)_{\alpha} \xrightarrow{\varphi\circ-} \Hom^{\ast}_{A_{\pal}}(\sQ,\Upsilon^{\ast}\sI)_{\alpha} \]
for every $\alpha\in\sN$.
\end{remark}

\begin{theorem}\label{thm.REndpseudomodules}
Let $\sF$ be a quasi-coherent sheaf on $X$, and let $\varepsilon\colon\sQ\to \Upsilon^{\ast}\sF$ be a cofibrant replacement in $\PM(A_{\pal})$. Then $\R\End(\sF)$ is represented by $\End_{A_{\pal}}^{\ast}(\sQ)$.
\begin{proof}
First notice that $\Hom^{\ast}_{A_{\pal}}(\sQ,\Upsilon^{\ast}\sI) \cong\Hom^{\ast}_{\Oh_X}(\Upsilon_!\sQ,\sI)$, the proof being the same as Proposition~\ref{prop.shriekadjoint}.
Now the statement follows by Proposition~\ref{proposition.REnd2} and Proposition~\ref{prop.REnd}.
\end{proof}
\end{theorem}

\bigskip

\subsection{$\R\End(\sF)$ via Thom-Whitney totalization}\label{section.TW}
The aim of this subsection is to prove that given a cofibrant replacement $\sQ\to \Upsilon^{\ast}\sF$ in $\PM(A_{\pal})$, then the derived endomorphisms of $\sF$ are represented by the Thom-Whitney totalization of a certain semicosimplicial DG-Lie algebra described in terms of $\sQ$, see Definition~\ref{def.semicosimplicialL}.

We begin by recalling the following construction. Let $\{U_j\}_{j\in J}$ be an affine open covering for a finite-dimensional Noetherian separated scheme $X$. Define
\[ \bar{\sN}_n=\{ (j_0,\dots,j_n) \in J^n \vert U_{j_0}\cap\dots\cap U_{j_n} \neq \emptyset \} \]
for any $n\in\N$. The \textbf{ordered nerve} of $\{U_j\}$ is the disjoint union $\bar{\sN} = \coprod\limits_{n\geq 0}\bar{I}_n$.
Notice that there exists a map
\[ \bar{\sN}\to \sN \, , \qquad \qquad \bar{\alpha}=(j_0,\dots,j_n)\mapsto \alpha=\{j_0,\dots,j_n\} \]
where $\sN$ is the nerve of $\{U_j\}$.

Consider $\sQ\in\PM(A_{\pal})$, and for every $n\in\N$ define
\[ \mathfrak{L}_n = \prod\limits_{\bar{\alpha}\in\bar{\sN}n}\Hom^{\ast}_{A_{\pal}}(\sQ,\sQ)_{\alpha} \]
where the product is taken in the category of DG-vector spaces.
Notice that $\mathfrak{L}^n$ is a DG-Lie algebra since every $\Hom^{\ast}_{A_{\pal}}(\sQ,\sQ)_{\alpha}\subseteq\prod\limits_{\gamma\geq\alpha}\Hom^{\ast}_{A_{\gamma}}(\sQ_{\gamma},\sQ_{\gamma})$ inherits a DG-Lie structure, where the bracket is the (graded) commutator. Moreover, for every monotone map $f\colon[n]\to[m]$ it is induced a map
\[ h_f\colon \bar{\sN}_m\to\bar{\sN}_n \, , \qquad \qquad \bar{\alpha}=(a_0,\dots,a_m) \mapsto h_f(\bar{\alpha})=(a_{f(0)},\dots,a_{f(n)}) \]
satisfying $h_f(\bar{\alpha})\leq\bar{\alpha}$ for every $\bar{\alpha}\in\bar{\sN}$. This in turn gives a map
\[ f_{\ast}= \left\{ f_{\bar{\beta}} \right\}_{\bar{\beta}\in \bar{\sN}_{m}} \colon \mathfrak{L}_n\to \mathfrak{L}_m \quad \mbox{defined by} \quad f_{\bar{\beta}}\left( \{\varphi_{\bar{\alpha}}\}_{\bar{\alpha}\in\bar{\sN}_n} \right) = \pi_{h_f(\bar{\beta})\bar{\beta}}\left( \varphi_{h_f(\bar{\beta})} \right) \in \mathfrak{L}_n \, , \]
where $\pi_{h_f(\bar{\beta})\bar{\beta}}\colon \Hom^{\ast}_{A_{\pal}}(\sQ,\sQ)_{h_f(\bar{\beta})}\to \Hom^{\ast}_{A_{\pal}}(\sQ,\sQ)_{\beta}$ is the natural projection.

\begin{definition}\label{def.semicosimplicialL}
For every $n\in\N$ and every $0\leq k\leq n+1$, define $\delta^k\colon [n]\to[n+1]$ as
\[ \delta^k(p)= \begin{cases}
p & \mbox{if } p<k \\
p+1 & \mbox{if } p\geq k
\end{cases} \]
Then the maps $\delta^k_{\ast}$ induce the semicosimplicial DG-Lie algebra
\[ \mathfrak{L} \colon \qquad \xymatrix{ \mathfrak{L}_0 \ar@<0.5ex>[r]\ar@<-0.5ex>[r] & \mathfrak{L}_1 \ar@<0.9ex>[r]\ar[r]\ar@<-0.9ex>[r] & \mathfrak{L}_1 \ar@<1.2ex>[r]\ar@<0.4ex>[r]\ar@<-0.4ex>[r]\ar@<-1.2ex>[r] & \cdots   } \]
\end{definition}

Similarly we now introduce three semicosimplicial complexes. Let $\sQ\to \Upsilon^{\ast}\sF$ be a cofibrant replacement for $\Upsilon^{\ast}\sF$ in $\PM(A_{\pal})$ and consider an injective resolution $\sF\to \sI$, then define
\[ \mathfrak{B}^{\sQ\sF}  \colon \qquad \xymatrix{ \mathfrak{B}_0^{\sQ\sF} = \prod\limits_{\bar{\alpha}\in\bar{\sN}_0}\Hom_{A_{\pal}}^{\ast}(\sQ,\Upsilon^{\ast}\sF)_{\alpha} \ar@<0.5ex>[r]\ar@<-0.5ex>[r] & \mathfrak{B}_1^{\sQ\sF} = \prod\limits_{\bar{\alpha}\in\bar{\sN}_1}\Hom_{A_{\pal}}^{\ast}(\sQ,\Upsilon^{\ast}\sF)_{\alpha} \ar@<0.9ex>[r]\ar[r]\ar@<-0.9ex>[r] & \cdots   } \]
\[ \mathfrak{B}^{\sQ\sI}  \colon \qquad \xymatrix{ \mathfrak{B}_0^{\sQ\sI} = \prod\limits_{\bar{\alpha}\in\bar{\sN}_0}\Hom_{A_{\pal}}^{\ast}(\sQ,\Upsilon^{\ast}\sI)_{\alpha} \ar@<0.5ex>[r]\ar@<-0.5ex>[r] & \mathfrak{B}_1^{\sQ\sI} = \prod\limits_{\bar{\alpha}\in\bar{\sN}_1}\Hom_{A_{\pal}}^{\ast}(\sQ,\Upsilon^{\ast}\sI)_{\alpha} \ar@<0.9ex>[r]\ar[r]\ar@<-0.9ex>[r] & \cdots   } \]
\[ \mathfrak{B}^{\sF\sI} \colon \qquad \xymatrix{ \mathfrak{B}_0^{\sF\sI} =\prod\limits_{\bar{\alpha}\in\bar{\sN}_0}\Hom^{\ast}_{\Oh_X}(i_{\alpha!}(\sF\vert_{U_{\alpha}}),\sI) \ar@<0.5ex>[r]\ar@<-0.5ex>[r] & \mathfrak{B}_1^{\sF\sI}=\prod\limits_{\bar{\alpha}\in\bar{\sN}_1}\Hom^{\ast}_{\Oh_X}(i_{\alpha!}(\sF\vert_{U_{\alpha}}),\sI) \ar@<0.9ex>[r]\ar[r]\ar@<-0.9ex>[r] & \cdots   } \]
where we denoted by $i_{\alpha}\colon U_{\alpha}\to X$ the natural inclusion.
Notice that the maps defined in Proposition~\ref{prop.REnd} and in Proposition~\ref{proposition.REnd2} induce semicosimplicial morphisms
\[ \mathfrak{L} \to \mathfrak{B}^{\sQ\sF} \to \mathfrak{B}^{\sQ\sI} \leftarrow \mathfrak{B}^{\sF\sI} \; . \]

Recall that for a semicosimplicial DG-vector space $V$ the Thom-Whitney-Sullivan totalization is the DG-vector space defined by
\[ \Tot_{TW}(V) = \left\{ (x_n)\in \prod\limits_{n\geq 0} \Omega_n\otimes V_n \,\Bigl\vert\, (\delta_k^{\ast}\otimes\Id)x_n = (\Id\otimes\delta_k)x_{n-1} \mbox{ for every } 0\leq k\leq n\right\} \]
where $\Omega_n=\frac{\K[t_0,\dots,t_n,dt_0,\dots,dt_n]}{(\sum t_i-1,\sum dt_i)}$ is the graded algebra of polynomial differential forms on the $n$-simplex.
Moreover, to every semicosimplicial DG-vector space $V$ is associated the complex
\[ C(V) = \bigoplus\limits_{p\in\N}\prod\limits_{n\in\N} V_n[-n]^p = \bigoplus\limits_{p\in\N}\prod\limits_{n\in\N}V_n^{p-n} \]
which is quasi-isomorphic to the totalization via the Whitney integration map $\int\colon \Tot_{TW}(V) \to C(V)$, see~\cite{Whi}. Given a map of DG-vector spaces $g\colon W\to V_0$ satisfying $\delta_0g=\delta_1g$, it is induced a morphism $\hat{g}\colon W\to \Tot_{TW}(V)$ defined by $\hat{g}(w)=(1\otimes g(w),1\otimes \delta_0g(w),1\otimes \delta_0^2g(w),\dots)$. Using the semicosimplicial identities it is straightforward to prove that the composition $\int\circ g$ is in fact the composition of $g$ with the natural inclusion $V_0\to C(V)$. In this way it is induced a natural map
\[ \Hom^{\ast}_{A_{\pal}}(\sQ,\sQ) \to \Tot_{TW}(\mathfrak{L}) \]
which respects the DG-Lie structure.

The aim of this subsection is to prove that $\Hom^{\ast}_{A_{\pal}}(\sQ,\sQ) \to \Tot_{TW}(\mathfrak{L})$ is a quasi-isomorphism of DG-Lie algebras. Actually we shall prove much more: there exists a commutative diagram
\begin{equation}\label{equation.diagram}
\xymatrix{	\Hom^{\ast}_{A_{\pal}}(\sQ,\sQ) \ar@{->}[r] \ar@{->}[d] & \Hom^{\ast}_{A_{\pal}}(\sQ,\Upsilon^{\ast}\sF) \ar@{->}[r] \ar@{->}[d] & \Hom^{\ast}_{A_{\pal}}(\sQ,\Upsilon^{\ast}\sI) \ar@{->}[d] & \Hom^{\ast}_{\Oh_X}(\sF,\sI) \ar@{->}[l] \ar@{->}[d] \\
\Tot_{TW}(\mathfrak{L})  \ar@{->}[r] \ar@{->}[d]^{\int} & \Tot_{TW}(\mathfrak{B}^{\sQ\sF})  \ar@{->}[r] \ar@{->}[d]^{\int} & \Tot_{TW}(\mathfrak{B}^{\sQ\sI}) \ar@{->}[d]^{\int} & \Tot_{TW}(\mathfrak{B}^{\sF\sI})  \ar@{->}[l] \ar@{->}[d]^{\int} \\
C(\mathfrak{L})  \ar@{->}[r] & C(\mathfrak{B}^{\sQ\sF})  \ar@{->}[r] & C(\mathfrak{B}^{\sQ\sI}) & C(\mathfrak{B}^{\sF\sI})  \ar@{->}[l] 	\; .}
\end{equation}
where all maps are quasi-isomorphisms.

\begin{lemma}\label{lemma.verticalquasiiso}
The vertical map $\Hom^{\ast}_{\Oh_X}(\sF,\sI)\xrightarrow{\xi} \Tot_{TW}(\mathfrak{B}^{\sF\sI})$ appearing in diagram~\eqref{equation.diagram} is a quasi-isomorphism.
\begin{proof}
As already noticed above the Whitney integration map $\int\colon\Tot_{TW}(\mathfrak{B}^{\sF\sI})\to C(\mathfrak{B}^{\sF\sI})$ is a quasi-isomorphism. Therefore, in order to prove the statement it is sufficient to show that the composition $\int\circ\xi$ is an isomorphism in cohomology.
To this aim we introduce two double complexes
\[ A_{ij} =\begin{cases}
\Hom^i_{\Oh_X}(\sF,\sI) & \mbox{ if } j=0 \\
0 & \mbox{ otherwise }
\end{cases} \qquad
B_{ij} =\prod\limits_{\bar{\alpha}\in\bar{\sN}_j}\Hom_{\Oh_{U_{\alpha}}}^i\left(\sF\vert_{U_{\alpha}},\sI\vert_{U_{\alpha}}\right) \]
defined for $i,j\geq 0$. Restrictions give a map of double complexes $\{A_{ij}\to B_{ij}\}_{i,j\geq 0}$, which in turn corresponds to a morphism between the associated complexes
\[ f\colon A^{\cdot}=\bigoplus\limits_{n\in\N} \Hom^n_{\Oh_X}(\sF,\sI) \to B^{\cdot} = \bigoplus\limits_{n\in\N}\bigoplus\limits_{i=0}^{n} B_{n-i,i} \; .\]
Now, consider the following complete and exhaustive filtrations
\[ F^pA^{\cdot} = \bigoplus\limits_{i\geq p}\Hom^i_{\Oh_X}(\sF,\sI) \; , \qquad \qquad F^pB^{\cdot}=\bigoplus\limits_{i\geq p}\bigoplus\limits_{j\geq 0}B_{i,j} \; , \qquad \qquad p\in\N \]
together with the induced morphism
\[ \hat{f}^p\colon \frac{F^pA^{\cdot}}{F^{p+1}A^{\cdot}} = \Hom_{\Oh_X}^p(\sF,\sI) \longrightarrow \bigoplus\limits_{j\geq 0}\prod\limits_{\bar{\alpha}\in\bar{\sN}_j}\Hom_{\Oh_{U_{\alpha}}}^p\left(\sF\vert_{U_{\alpha}},\sI\vert_{U_{\alpha}}\right) = \frac{F^pB^{\cdot}}{F^{p+1}B^{\cdot}} \; . \]
Observe that for every $p\in\N$ the map $\hat{f}^p$ is a quasi-isomorphism; in fact by the degreewise injectivity of $\sI$ it follows that the restriction map
\[ \Hom^p_{\Oh_X}(\sF,\sI) \to \Hom_{\Oh_X}^p(i_!\sF,\sI)=\Hom_{\Oh_X\vert_U}^p\left(\sF\vert_V,\sI\vert_V\right) \]
is surjective for every open subset $i\colon V\to X$, therefore the sequence
\[ 0\to \Hom_{\Oh_X}^p(\sF,\sI) \to \prod\limits_{\bar{\alpha}\in\bar{\sN}_0}\Hom^p_{\Oh_{U_{\alpha}}}\left(\sF\vert_{U_{\alpha}},\sI\vert_{U_{\alpha}}\right) \to \prod\limits_{\bar{\beta}\in\bar{\sN}_1}\Hom^p_{\Oh_{U_{\beta}}}\left(\sF\vert_{U_{\beta}},\sI\vert_{U_{\beta}}\right) \to \cdots \]
is exact because flasque sheaves are acyclic.
It follows that the map $f\colon A^{\cdot}\to B^{\cdot}$ is a quasi-isomorphism.

To conclude the proof it is sufficient to observe that $f$ is indeed the composition $\int\circ\xi$. Clearly $A^{\cdot}=\Hom_{\Oh_X}^{\ast}(\sF,\sI)$; moreover
\[ B^{\cdot} = \bigoplus\limits_{n\in\N}\bigoplus\limits_{i=0}^{n} B_{n-i,i} = \bigoplus\limits_{n\in\N}\bigoplus\limits_{i=0}^{n} \prod\limits_{\bar{\alpha}\in\bar{\sN}_i}\Hom_{\Oh_{U_{\alpha}}}^{n-i}\left(\sF\vert_{U_{\alpha}},\sI\vert_{U_{\alpha}}\right) =  \bigoplus\limits_{n\in\N}\prod\limits_{i=0}^{n} \Hom_{\Oh_{U_{\alpha}}}^{n-i}\left(\bigoplus\limits_{\bar{\alpha}\in\bar{\sN}_i}i_{\alpha!}(\sF\vert_{U_{\alpha}}),\sI\right) \]
so that $B^{\cdot}=C(\mathfrak{B}^{\sF\sI})$. Now, the map $\int\circ\xi$ is the same as the composition
\[ \Hom_{\Oh_X}^{\ast}(\sF,\sI) \to \prod\limits_{\alpha\in\sN_0}\Hom_{\Oh_{U_{\alpha}}}^{\ast}\left(\sF\vert_{U_{\alpha}},\sI\vert_{U_{\alpha}}\right) \to C(\mathfrak{B}^{\sF\sI}) \]
which is precisely $f$ as claimed.
\end{proof}
\end{lemma}

\begin{theorem}\label{thm.diagram}
All the maps appearing in diagram~\eqref{equation.diagram} are quasi-isomorphisms.
\begin{proof}
The maps in the first row have been discussed in Proposition~\ref{prop.REnd} and Proposition~\ref{rmk.REnd2}.
Now, recall that to prove that the map between complexes associated to semicosimplicial DG-vector spaces is a quasi-isomorphisms, it is sufficient to prove that it is induced by a semicosimplicial quasi-isomorphism between them. By Remark~\ref{rmk.REnd2} and by Proposition~\ref{prop.REnd} there are quasi-isomorphisms
\[ \xymatrix{	\Hom^{\ast}_{A_{\pal}}(\sQ,\sQ)_{\alpha} \ar@{->}[r]^{\varepsilon\circ-} & \Hom^{\ast}_{A_{\pal}}(\sQ,\Upsilon^{\ast}\sF)_{\alpha} \ar@{->}[r]^{\varphi\circ-} & \Hom^{\ast}_{A_{\pal}}(\sQ,\Upsilon^{\ast}\sI)_{\alpha} \ar@{->}[d]^{\cong} \\
 & & \Hom^{\ast}_{\Oh_{U_{\alpha}}}(\Upsilon_!\sQ,\sI) & \Hom^{\ast}_{\Oh_{U_{\alpha}}}\left(\sF\vert_{U_{\alpha}},\sI\vert_{U_{\alpha}}\right) \ar@{->}[l]_{-\circ\varepsilon} 	} \]
for every $\alpha\in\sN$, which induce semicosimplicial quasi-isomorphisms $\mathfrak{L} \to \mathfrak{B}^{\sQ\sF}  \to \mathfrak{B}^{\sQ\sI} \leftarrow \mathfrak{B}^{\sF\sI}$.
Therefore the maps in the bottom row are all quasi-isomorphisms. Moreover, since for every DG-vector space $V$ the map $\int\colon\Tot_{TW}(V)\to C(V)$ is a quasi-isomorphism, by the \emph{2 out of 3} property also the maps in the middle row are quasi-isomorphisms.

To conclude the proof recall that the map $\Hom^{\ast}_{\Oh_X}(\sF,\sI)\to \Tot_{TW}(\mathfrak{B}^{\sF\sI})$ is a quasi-isomorphism by Lemma~\ref{lemma.verticalquasiiso}. Hence the statement follows again by the \emph{2 out of 3} property.
\end{proof}
\end{theorem}

\begin{corollary}\label{corollary.diagram}
Let $\sF$ be a quasi-coherent sheaf on $X$, and let $\varepsilon\colon\sQ\to \Upsilon^{\ast}\sF$ be a cofibrant replacement in $\PM(A_{\pal})$. Denote by $\mathfrak{L}$ the semicosimplicial DG-Lie algebra introduced in Definition~\ref{def.semicosimplicialL}. Then $\R\End(\sF)$ is represented by $\Tot_{TW}(\mathfrak{L})$.
\begin{proof}
Immediate consequence of Theorem~\ref{thm.REndpseudomodules} and Theorem~\ref{thm.diagram}
\end{proof}
\end{corollary}

\begin{remark}\label{rmk.diagram}
Another consequence of Theorem~\ref{thm.diagram} is the existence of a quasi-isomorphism of differential graded Lie algebras $\Hom^{\ast}_{A_{\pal}}(\sQ,\sQ) \to \Tot_{TW}(\mathfrak{L})$. This implies that the associated deformations functors defined through Maurer-Cartan elements modulo gauge equivalence are isomorphic:
\[ \Def_{\Hom^{\ast}_{A_{\pal}}(\sQ,\sQ)} \cong\Def_{\Tot_{TW}(\mathfrak{L})} \]
see~\cite[Corollary 5.52]{Man3}.
\end{remark}

\bigskip

\subsection{$\R\End(\sF)$ in presence of a locally free resolution}\label{section.locallyfree}

Let $\sE\to\sF$ be a locally free resolution for a quasi-coherent sheaf $\sF$ over $X$. Recall that if $X$ is smooth projective such a resolution always exists, but we keep working in full generality only assuming $X$ to be a finite-dimensional separated Noetherian scheme over $\K$. Moreover we choose an affine open cover $\{U_h\}_{h\in H}$ for $X$ such that the restriction $\sE\vert_{U_{\alpha}}$ is a complex of free sheaves for every $\alpha\in\sN$. Notice that:
\begin{enumerate}
\item $\Upsilon^{\ast}\sE\in\PM(A_{\pal})$ is quasi-coherent,
\item $(\Upsilon^{\ast}\sE)_{\alpha}$ is cofibrant in $\DGMod(A_{\alpha})$ for every $\alpha\in\sN$,
\item $\Upsilon^{\ast}\sE$ is not necessarily cofibrant in $\PM(A_{\pal})$.
\end{enumerate}

\begin{lemma}\label{lemma.locallyfree}
Let $\sE\to\sF$ be a locally free resolution, and consider a cofibrant replacement $\sQ\xrightarrow{\pi} \Upsilon^{\ast}\sE$ in $\PM(A_{\pal})$. Fix $\alpha\in\sN$; then all the maps in the commutative square
\[ \xymatrix{		\Hom_{A_{\pal}}^{\ast}(\sQ,\Upsilon^{\ast}\sF)_{\alpha} \ar@{->}[d] & \Hom_{A_{\pal}}^{\ast}(\Upsilon^{\ast}\sE,\Upsilon^{\ast}\sF)_{\alpha} \ar@{->}[l]\ar@{->}[d] \\
\Hom^{\ast}_{A_{\alpha}}(\sQ_{\alpha},(\Upsilon^{\ast}\sF)_{\alpha}) & \Hom^{\ast}_{A_{\alpha}}\left((\Upsilon^{\ast}\sE)_{\alpha},(\Upsilon^{\ast}\sF)_{\alpha}\right) \ar@{->}[l] 	} \]
are quasi-isomorphisms, where the vertical arrows are the natural projections.
\begin{proof}
First notice that the vertical arrow on the right is clearly an isomorphism.
Moreover, the bottom arrow is a quasi-isomorphism because it is induced by the map $\sQ_{\alpha}\to(\Upsilon^{\ast}\sE)_{\alpha}$, which is a weak equivalence between cofibrant objects in $\DGMod(A_{\alpha})$. By the \emph{2 out of 3} axiom it is then sufficient to prove that the projection
\[ \pi\colon\Hom_{A_{\pal}}^{\ast}(\sQ,\Upsilon^{\ast}\sF)_{\alpha}\to \Hom^{\ast}_{A_{\alpha}}(\sQ_{\alpha},(\Upsilon^{\ast}\sF)_{\alpha}) \]
is a quasi-isomorphism. We begin by showing the surjectivity in cohomology. To this aim, take $\varphi_{\alpha}\in Z^0\left( \Hom^{\ast}_{A_{\alpha}}(\sQ_{\alpha},(\Upsilon^{\ast}\sF)_{\alpha}) \right) = \Hom_{A_{\alpha}}(\sQ_{\alpha},(\Upsilon^{\ast}\sF)_{\alpha})$. By induction, fix $\beta\in\sN$ such that $\alpha<\beta$ and suppose we have already constructed maps $\varphi_{\gamma}\in\Hom_{A_{\gamma}}(\sQ_{\gamma},(\Upsilon^{\ast}\sF)_{\gamma})$ for every
\[ \gamma\in \sR_{\alpha\beta}=\{\gamma\in\sN\,\vert\,\alpha\leq\gamma<\alpha\} \]
satisfying the necessary commutativity relations. In order to define $\varphi_{\beta}\in\Hom_{A_{\beta}}(\sQ_{\beta},(\Upsilon^{\ast}\sF)_{\beta})$ first notice that the map
\[ \colim_{\gamma\in\sR_{\alpha\beta}}(\sQ_{\gamma}\otimes_{A_{\gamma}}A_{\beta}) \to \sQ_{\beta} \]
is a cofibration in $\DGMod(A_{\beta})$ by Remark~\ref{rmk.restrictioncofibrant}. Notice that $\sQ$ is a quasi-coherent $A_{\pal}$-module by Remark~\ref{rmk.qcoh-homotopyinvariance}, so that the map
\[ \left\{ \sQ_{\gamma}\otimes_{A_{\gamma}}A_{\beta} \to \sQ_{\beta} \right\}_{\gamma\in\sR_{\alpha\beta}} \]
is a Reedy weak equivalence. Moreover, the diagram $\left\{ \sQ_{\gamma}\otimes_{A_{\gamma}}A_{\beta} \right\}_{\gamma\in\sR_{\alpha\beta}}$ is Reedy cofibrant by Remark~\ref{rmk.restrictioncofibrant}, and $\left\{ \sQ_{\beta} \right\}_{\gamma\in\sR_{\alpha\beta}}$ is Reedy cofibrant since $\sR_{\alpha\beta}$ is connected. It follows that the map
\[ \colim_{\gamma\in\sR_{\alpha\beta}}(\sQ_{\gamma}\otimes_{A_{\gamma}}A_{\beta}) \to \colim_{\gamma\in\sR_{\alpha\beta}}\sQ_{\beta} \cong \sQ_{\beta} \]
is a weak equivalence since the left Quillen functor $\colim\colon\DGMod(A_{\beta})^{\sR_{\alpha\beta}}\to \DGMod(A_{\beta})$ preserves weak equivalences between Reedy cofibrant objects by Ken Brown's Lemma. Hence the diagram
\[ \xymatrix{	\colim\limits_{\gamma\in\sR_{\alpha\beta}}(\sQ_{\gamma}\otimes_{A_{\gamma}}A_{\beta}) \ar@{->}[d]_{\sC\sW} \ar@{->}[r] & (\Upsilon^{\ast}\sF)_{\beta} \\
\sQ_{\beta} \ar@{.>}[ur]_{\varphi_{\beta}} & 		} \]
admits the required dotted lifting. This proves that $\pi$ is surjective in cohomology in degree $0$. For the general case it is sufficient to observe that
\[ Z^n\left(\Hom^{\ast}_{A_{\alpha}}(\sQ_{\alpha},(\Upsilon^{\ast}\sF)_{\alpha})\right) \cong Z^0\left( \Hom^{\ast}_{A_{\alpha}}(\sQ_{\alpha},(\Upsilon^{\ast}\sF)_{\alpha}[n]) \right) \, .\]
We are left with the proof of the injectivity of $\pi$ in cohomology. To this aim, take $\{\varphi_{\gamma}\}_{\gamma\geq\alpha}$ in $\Hom_{A_{\pal}}(\sQ,\Upsilon^{\ast}\sF)_{\alpha}$ and suppose that $\varphi_{\alpha}\colon\sQ_{\alpha}\to (\Upsilon^{\ast}\sF)_{\alpha}$ is homotopic to the zero map; i.e. $\pi(\{\varphi_{\gamma}\})=0$ in $H^0\left( \Hom^{\ast}_{A_{\alpha}}(\sQ_{\alpha},(\Upsilon^{\ast}\sF)_{\alpha})\right)$. This is equivalent to say that the diagram of solid arrows
\[ \xymatrix{	 & \cone\left( \Id_{(\Upsilon^{\ast}\sF)_{\alpha}[-1]} \right) \ar@{->}[d]^{p_{\alpha}} \\
\sQ_{\alpha} \ar@{->}[r]_{\varphi_{\alpha}} \ar@{.>}[ur]^{h_{\alpha}} & (\Upsilon^{\ast}\sF)_{\alpha} 	} \]
admits the dotted lifting $h_{\alpha}$. Recall that
\[ \cone\left( \Id_{(\Upsilon^{\ast}\sF)_{\alpha}[-1]} \right) = (\Upsilon^{\ast}\sF)_{\alpha}\oplus(\Upsilon^{\ast}\sF)_{\alpha}[-1] \]
as graded $A_{\alpha}$-modules, and $p_{\alpha}$ is the projection on the first summand. In order to prove that $\{\varphi_{\gamma}\}$ is exact we proceed by induction: fix $\beta\in\sN$ such that $\alpha<\beta$ and suppose that the homotopy $h_{\alpha}$ has been lifted to $h_{\gamma}\colon\sQ_{\gamma}\to \cone\left(\Id_{(\Upsilon^{\ast}\sF)_{\gamma}[-1]}\right)$ for every $\gamma\in\sR_{\alpha\beta}=\{\gamma\in\sN\,\vert\,\alpha\leq\gamma<\beta\}$. We need to prove the existence of the dotted lifting in the diagram below
\[ \xymatrix{	 & & \colim\limits_{\gamma\in\sR_{\alpha\beta}} \cone\left(\Id_{(\Upsilon^{\ast}\sF)_{\gamma}[-1]}\right) \ar@{->}[d] \ar@{->}[drr] & & \\
\colim\limits_{\gamma\in\sR_{\alpha\beta}}\sQ_{\gamma}\otimes_{A_{\gamma}}A_{\beta} \ar@{->}[drr] \ar@{->}[urr]^{\hat{h}} \ar@{->}[rr]^{\hat{\varphi}} & & (\Upsilon^{\ast}\sF)_{\beta} \ar@{->}[drr]|-{ } & & \cone\left(\Id_{(\Upsilon^{\ast}\sF)_{\beta}[-1]}\right) \ar@{->}[d]^{p_{\beta}} \\
 & & \sQ_{\beta} \ar@{.>}[urr]^<<<<<<<<{h_{\beta}} \ar@{->}[rr]_{\varphi_{\beta}} & & (\Upsilon^{\ast}\sF)_{\beta}	} \]
where $\hat{h}$ is induced by $\{h_{\gamma}\}_{\gamma\in\sR_{\alpha\beta}}$ and $\hat{\varphi}$ is induced by $\{\varphi_{\gamma}\}_{\gamma\in\sR_{\alpha\beta}}$. Notice that $p_{\beta}$ is surjective (hence a fibration), and $\colim\limits_{\gamma\in\sR_{\alpha\beta}}\sQ_{\gamma}\otimes_{A_{\gamma}}A_{\beta} \to \sQ_{\beta}$ is a trivial cofibration as proved above; therefore the statement follows by the lifting property.
\end{proof}
\end{lemma}

\begin{remark}\label{rmk.projectionqiso}
Even if $\sF$ does not admit a locally free resolution, we can consider a cofibrant replacement $\sQ\to \Upsilon^{\ast}\sF$ in $\PM(A_{\pal})$: the same argument of Lemma~\ref{lemma.locallyfree} shows that the projection
\[ \Hom_{A_{\pal}}^{\ast}(\sQ,\Upsilon^{\ast}\sF)_{\alpha} \to \Hom^{\ast}_{A_{\alpha}}(\sQ_{\alpha},(\Upsilon^{\ast}\sF)_{\alpha}) \]
is a quasi-isomorphism.
\end{remark}

\begin{remark}\label{rmk.qisolocallyfree}
In the proof of Lemma~\ref{lemma.locallyfree}, the fact that $\Upsilon^{\ast}\sF$ is concentrated in degree $0$ does not play any role. Therefore for every $\alpha\in\sN$ the same argument leads to a quasi-isomorphism
\[ -\circ\pi\colon \Hom^{\ast}_{A_{\pal}}(\Upsilon^{\ast}\sE,\Upsilon^{\ast}\sE)_{\alpha} \to \Hom^{\ast}_{A_{\pal}}(\sQ,\Upsilon^{\ast}\sE)_{\alpha} \]
where $\pi\colon\sQ\to\Upsilon^{\ast}\sE$ is a cofibrant replacement in $\PM(A_{\pal})$.
\end{remark}

Given a locally free resolution $\sE\to\sF$ on $X$, we consider the associated \v{C}ech semicosimplicial DG-Lie algebra
\[ \mathfrak{h} \colon \qquad \xymatrix{ \prod\limits_{\bar{\alpha}\in\bar{\sN}_0}\Hom^{\ast}_{\Oh_{U_{\alpha}}}(\sE\vert_{U_{\alpha}},\sE\vert_{U_{\alpha}}) \ar@<0.5ex>[r]\ar@<-0.5ex>[r] & \prod\limits_{\bar{\beta}\in\bar{\sN}_1}\Hom^{\ast}_{\Oh_{U_{\beta}}}(\sE\vert_{U_{\beta}},\sE\vert_{U_{\beta}}) \ar@<0.9ex>[r]\ar[r]\ar@<-0.9ex>[r] & \cdots   } \]
which will give us another model for derived endomorphisms of $\sF$.

\begin{theorem}\label{thm.locallyfreeREnd}
Let $\sF$ be a quasi-coherent sheaf on $X$, and let $\sE\to\sF$ be a locally free resolution. Denote by $\mathfrak{h}$ the \v{C}ech semicosimplicial DG-Lie algebra as above. Then $\R\End(\sF)$ is represented by $\Tot_{TW}(\mathfrak{h})$.
\begin{proof}
Take a cofibrant replacement $\sQ\to\Upsilon^{\ast}\sE$ in $\PM(A_{\pal})$ and fix $\alpha\in\sN$. By Lemma~\ref{lemma.locallyfree} there exists a quasi-isomorphism
\[ \Hom_{A_{\pal}}^{\ast}(\sQ,\Upsilon^{\ast}\sF)_{\alpha} \leftarrow \Hom_{A_{\pal}}^{\ast}(\Upsilon^{\ast}\sE,\Upsilon^{\ast}\sF)_{\alpha} \cong \Hom^{\ast}_{A_{\alpha}}\left((\Upsilon^{\ast}\sE)_{\alpha},(\Upsilon^{\ast}\sF)_{\alpha}\right) \cong \Hom^{\ast}_{\Oh_{U_{\alpha}}}(\sE\vert_{U_{\alpha}},\sF\vert_{U_{\alpha}}) \; . \]
Moreover the map $\Hom_{\Oh_{U_{\alpha}}}^{\ast}(\sE\vert_{U_{\alpha}},\sE\vert_{U_{\alpha}})\to \Hom_{\Oh_{U_{\alpha}}}^{\ast}(\sE\vert_{U_{\alpha}},\sF\vert_{U_{\alpha}})$ is a quasi-isomorphism, being $\sE\vert_{U_{\alpha}}$ a complex of free sheaves. Therefore we obtain a quasi-isomorphism
\[ \Hom_{\Oh_{U_{\alpha}}}^{\ast}(\sE\vert_{U_{\alpha}},\sE\vert_{U_{\alpha}}) \to \Hom_{A_{\pal}}^{\ast}(\sQ,\Upsilon^{\ast}\sF)_{\alpha} \]
which extends to a semicosimplicial quasi-isomorphism $\mathfrak{h}\to\mathfrak{B}^{\sQ\sF}$, so that the induced map $\Tot_{TW}(\mathfrak{h})\to\Tot_{TW}(\mathfrak{B}^{\sQ\sF})$ is a quasi-isomorphism. The statement follows by Theorem~\ref{thm.diagram} and Corollary~\ref{corollary.diagram}.
\end{proof}
\end{theorem}

Theorem~\ref{rmk.qisolocallyfree} essentially states that $H^k\left(\Tot_{TW}(\mathfrak{h})\right) = \Ext^k_{\Oh_X}(\sF,\sF)$ for every $k\in\N$. For future purposes, we are now interested in a stronger result, namely that $\Tot_{TW}(\mathfrak{h})$, $\Tot_{TW}(\mathfrak{L})$ and $\End^{\ast}_{A_{\pal}}(\sQ)$ are quasi-isomorphic as DG-Lie algebras, so that in particular the associated deformation functors $\Def_{\Tot_{TW}(\mathfrak{h})}$, $\Def_{\Tot_{TW}(\mathfrak{L})}$ and $\Def_{\Hom^{\ast}_{A_{\pal}}(\sQ)}$ will be isomorphic to each other. Recall that it has been already proven in Subsection~\ref{section.TW} that $\Def_{\Tot_{TW}(\mathfrak{L})}\cong\Def_{\Hom^{\ast}_{A_{\pal}}(\sQ)}$.

\begin{lemma}\label{lemma.locallyfreeDGLA}
Let $\sE\to\sF$ be a locally free resolution, and consider a cofibrant replacement $\sQ\xrightarrow{\pi} \Upsilon^{\ast}\sE$ in $\PM(A_{\pal})$. Fix $\alpha\in\sN$ and define the DG-Lie algebra
\[ M_{\alpha}= \left\{ (f,g)\in\Hom^{\ast}_{A_{\pal}}(\sQ,\sQ)_{\alpha}\times\Hom^{\ast}_{A_{\pal}}(\Upsilon^{\ast}\sE,\Upsilon^{\ast}\sE)_{\alpha} \,\vert\, \pi\circ f=g\circ\pi \right\} \; . \]
Then there exists a commutative square
\[ \xymatrix{	M_{\alpha} \ar@{->}[rr]^-{p_2} \ar@{->}[d]_{p_1} & & \Hom^{\ast}_{A_{\pal}}(\Upsilon^{\ast}\sE,\Upsilon^{\ast}\sE)_{\alpha} \ar@{->}[d]^{-\circ\pi} \\
\Hom^{\ast}_{A_{\pal}}(\sQ,\sQ)_{\alpha} \ar@{->}[rr]_-{\pi\circ-} & & \Hom^{\ast}_{A_{\pal}}(\sQ,\Upsilon^{\ast}\sE)_{\alpha} 	} \]
where every map is a quasi-isomorphism.
\begin{proof}
First notice that the map
\[ \pi\circ-\colon \Hom^{\ast}_{A_{\pal}}(\sQ,\sQ)_{\alpha} \to \Hom^{\ast}_{A_{\pal}}(\sQ,\Upsilon^{\ast}\sE)_{\alpha} \]
is a quasi-isomorphism being $\sQ$ cofibrant in $\PM(A_{\pal})$, see Remark~\ref{rmk.restrictioncofibrant}. Moreover, the map
\[ -\circ\pi\colon \Hom^{\ast}_{A_{\pal}}(\Upsilon^{\ast}\sE,\Upsilon^{\ast}\sE)_{\alpha} \to\Hom^{\ast}_{A_{\pal}}(\sQ,\Upsilon^{\ast}\sE)_{\alpha} \]
is a quasi-isomorphism by Remark~\ref{rmk.qisolocallyfree}. By the functoriality of cohomology, to prove the statement it is sufficient to show that the projection $p_1$ is a quasi-isomorphism. To this aim, first observe that $\sQ$ is cofibrant and $\pi$ is surjective, so that the map $p_1$ is surjective by Lemma~\ref{lemma.cofibrantprojective}. Moreover, the complex $\ker(p_1) = \Hom^{\ast}_{A_{\pal}}(\sQ,\ker(\pi))_{\alpha}$ is acyclic, being $\sQ$ cofibrant and $\ker(\pi)$ acyclic. The statement follows.
\end{proof}
\end{lemma}

\begin{theorem}\label{thm.locallyfreeDef}
Let $\sE\to\sF$ be a locally free resolution, and consider a cofibrant replacement $\sQ\xrightarrow{\pi} \Upsilon^{\ast}\sE$ in $\PM(A_{\pal})$. Let $\mathfrak{L}$ be the semicosimplicial DG-Lie algebra associated to $\sQ$ as in Definition~\ref{def.semicosimplicialL}. Then $\Tot_{TW}(\mathfrak{L})$ and $\Tot_{TW}(\mathfrak{h})$ are quasi-isomorphic as DG-Lie algebras. In particular, the associated deformation functors $\Def_{\Tot_{TW}(\mathfrak{L})}$ and $\Def_{\Tot_{TW}(\mathfrak{h})}$ are naturally isomorphic.
\begin{proof}
It is sufficient to observe that by Lemma~\ref{lemma.locallyfreeDGLA} there exists quasi-isomorphisms
\[ \Hom^{\ast}_{A_{\pal}}(\sQ,\sQ)_{\alpha}\leftarrow M_{\alpha}\to \Hom^{\ast}_{A_{\pal}}(\Upsilon^{\ast}\sE,\Upsilon^{\ast}\sE)_{\alpha} \]
of DG-Lie algebras inducing quasi-isomorphisms of semicosimplicial DG-Lie algebras. To conclude the proof recall that the Whitney integration maps lift quasi-isomorphisms between complexes associated to semicosimplicial DG-Lie algebras to quasi-isomorphisms between their totalizations.
\end{proof}
\end{theorem}

\bigskip

\section{Infinitesimal deformations of quasi-coherent sheaves}\label{section.deformations}

It is well known that infinitesimal deformations of a coherent sheaf on a smooth projective variety are related to $\Ext^{\ast}(\sF,\sF)$, see e.g.~\cite{FIM}. Using results of Section~\ref{section.REnd}, our aim is now to prove that the DG-Lie algebras $\End^{\ast}_{A_{\pal}}(\sQ)=\Hom^{\ast}_{A_{\pal}}(\sQ,\sQ)$ and $\Tot_{TW}(\mathfrak{L})$ control infinitesimal deformations of a quasi-coherent sheaf $\sF$ over a finite-dimensional Noetherian separated scheme $X$. Here $\sQ\to\Upsilon^{\ast}\sF$ is any cofibrant replacement in $\PM(A_{\pal})$.

For the reader convenience, we briefly recall the definition of the deformation functor associated to infinitesimal deformations of $\sF$. A deformation of $\sF$ over $A\in\Art_{\K}$ is a morphism $\pi\colon\sF_A\to \sF$ of sheaves of $\Oh_X\otimes A$-modules over $X\times \Spec(A)$, with $\sF_A$ flat over $A$, such that the reduced map $\sF_A\otimes_A\K\to\sF$ is an isomorphism.
We say that two deformations $\sF_A$ and $\sF_A'$ are isomorphic if there exists an isomorphism of sheaves $\varphi\colon\sF_A\to \sF_A'$ such that $\pi'\circ\varphi=\pi$.
The functor of infinitesimal deformations of $\sF$ up to isomorphism is denoted by $\Def_{\sF}\colon\Art_{\K}\to\Set$.

The main result of this section will be the existence of natural isomorphisms
\[ \Def_{\sF}\cong\Def_{\Tot_{TW}(\mathfrak{L})}\cong\Def_{\End^{\ast}_{A_{\pal}}(\sQ)} \; . \]
We shall give different proofs.
First recall that by Remark~\ref{rmk.diagram} there exists a natural isomorphism $\Def_{\End^{\ast}_{A_{\pal}}(\sQ)}\to\Def_{\Tot_{TW}(\mathfrak{L})}$.
In Subsection~\ref{section.FIM} we will use a powerful result of~\cite{FIM}, which will lead us to a natural isomorphism $\Def_{\sF}\cong\Def_{\Tot_{TW}(\mathfrak{L})}$.
In Subsection~\ref{section.deformationpseudomod} we will give an explicit natural isomorphism $\Def_{\End^{\ast}_{A_{\pal}}(\sQ)}\to\Def_{\sF}$.

\bigskip

\subsection{Deformations via descent of Deligne groupoid}\label{section.FIM}

We begin by recalling the construction of the functors $Z^1_{\g}, \, H^1_{\g}\colon \Art_{\K}\to\Set$ for any given semicosimplicial DG-Lie algebra
\[ \g \colon \qquad \xymatrix{ \g_0 \ar@<0.5ex>[r]^{\partial_{0,1}}\ar@<-0.5ex>[r]_{\partial_{1,1}} & \g_1 \ar@<1.4ex>[rr]^{\partial_{0,2}}\ar[rr]|-{\partial_{1,2}}\ar@<-1.4ex>[rr]_{\partial_{2,2}} &  & \cdots   } \; . \]
For every $A\in\Art_{\K}$ define $Z^1_{\g}(A) \subseteq (\g_0^1\oplus\g_1^0)\otimes\mA$ to be the subset of elements $(l,m)\in (\g_0^1\oplus\g_1^0)\otimes\mA$ satisfying
\[ \begin{cases}
dl+\frac{1}{2}[l,l]=0 \\
\partial_{1,1}l = e^m\ast\partial_{0,1}l \\
\partial_{0,2}m\bullet(-\partial_{1,2}m)\bullet\partial_{2,2}m = dn+[\partial_{2,2}\partial_{0,1}l,n] & \text{ for some } n\in\g_2^{-1}\otimes\mA
\end{cases} \]
where $\ast$ denotes the gauge action and $\bullet$ denotes the Baker-Campbell-Hausdorff product; i.e. $x\bullet y=\log(e^xe^y)$. There is an equivalence relation on $Z^1_{\g}(A)$: two elements $(l_0,m_0),(l_1,m_1)\in Z^1_{\g}(A)$ are equivalent if and only if there exist $a\in\g_0^0\otimes\mA$ and $b\in\g_1^{-1}\otimes\mA$ such that
\[ \begin{cases}
e^a\ast l_0 = l_1 \\
-m_0\bullet (-\partial_{1,1}a)\bullet\partial_{0,1}a = db + [\partial_{0,1}l_0,b] \; .
\end{cases} \]
We shall denote by $\sim$ the equivalent relation defined above; the functor of Artin rings $H_{\g}^1\colon \Art_{\K}\to\Set$ is defined as $H_{\g}^1(A)=\faktor{Z_{\g}^1(A)}{\sim}$ for every $A\in\Art_{\K}$.
This functor extends the one defined in~\cite{FMM} for semicosimplicial Lie algebras.
It was proven in~\cite{FIM} that there exists a commutative diagram of functors
\[ \xymatrix{	\DGLA_{H^{\geq0}}^{\Delta} \ar@{->}[dr]_{H^1_{\cdot}} \ar@{->}[rr]^{\Tot_{TW}(\cdot)} & & \DGLA \ar@{->}[dl]^{\Def_{\cdot}} \\
 & \Set^{\Art_{\K}} & 		} \]
where $\DGLA_{H^{\geq0}}^{\Delta}$ is the category of semicosimplicial DG-Lie $\K$-algebras with cohomology concentrated in positive degrees, $\DGLA$ is the category of DG-Lie $\K$-algebras, and $\Set^{\Art_{\K}}$ is the category of functors $\Art_{\K}\to\Set$. Moreover, the functor $\Def_{\cdot}\colon\DGLA\to\Set^{\Art_{\K}}$ is defined by Maurer-Cartan solution modulo gauge equivalence.

Our strategy is now clear: we first need to show that the semicosimplicial DG-Lie algebra $\mathfrak{L}$ defined in~\ref{def.semicosimplicialL} has cohomology concentrated in positive degrees, i.e. $\mathfrak{L}\in\DGLA_{H^{\geq0}}^{\Delta}$, then we conclude by showing that $\Def_{\sF}\cong H^1_{\mathfrak{L}}$.

\begin{lemma}\label{lemma.positivecohomology}
Let $\sF$ be a quasi-coherent sheaf on $X$, and take a cofibrant replacement $\sQ\to\Upsilon^{\ast}\sF$ in $\PM(A_{\pal})$. Then the associated semicosimplicial DG-Lie algebra $\mathfrak{L}$ defined in~\ref{def.semicosimplicialL} belongs to $\DGLA_{H^{\geq0}}^{\Delta}$.
\begin{proof}
Fix $\alpha\in\sN$; we need to show that $\Hom_{A_{\pal}}^{\ast}(\sQ,\sQ)_{\alpha}$ is acyclic in negative degrees. Consider the composition
\[ \Hom_{A_{\pal}}^{\ast}(\sQ,\sQ) \to \Hom_{A_{\pal}}^{\ast}(\sQ,\Upsilon^{\ast}\sF) \to \Hom_{A_{\alpha}}^{\ast}\left(\sQ_{\alpha},\sF(U_{\alpha})\right) \] 
where the first map is a quasi-isomorphism by Proposition~\ref{proposition.REnd2}, and the second map is a quasi-isomorphism by Remark~\ref{rmk.projectionqiso}. Now consider a projective resolution $P^{\cdot}\to \sF(U_{\alpha})$, which in particular is a cofibrant replacement in $\DGMod(A_{\alpha})$, see e.g.~\cite[Lemma 2.3.6]{Hov99}. Therefore there exist a quasi-isomorphism $q\colon \sQ_{\alpha}\to P^{\cdot}$ lifting $\sQ_{\alpha}\to \sF(U_{\alpha})$. By Ken Brown's Lemma, the functor $\Hom_{A_{\alpha}}^{\ast}\left(-,\sF(U_{\alpha})\right)$ maps weak equivalences between cofibrant objects to quasi-isomorphisms, so that the induced map
\[ \Hom_{A_{\alpha}}^{\ast}\left(P^{\cdot},\sF(U_{\alpha})\right) \xrightarrow{-\circ q} \Hom_{A_{\alpha}}^{\ast}\left(\sQ_{\alpha},\sF(U_{\alpha})\right) \]
is a quasi-isomorphism. Now the statement follows since the complex $\Hom_{A_{\alpha}}^{\ast}\left(P^{\cdot},\sF(U_{\alpha})\right)$ does not have non-zero $n$-cocycles for $n<0$.
\end{proof}
\end{lemma}

Fix $\alpha\in\sN$ and $A\in\Art_{\K}$; a Maurer-Cartan element $\{l_{\beta}\}_{\beta\geq\alpha}\in\Hom_{A_{\pal}}^{1}(\sQ,\sQ)_{\alpha}\otimes\mA$ defines complexes $(\sQ_{\beta}\otimes A,d_{\sQ_{\beta}}+l_{\beta})$ for every $\beta\geq\alpha$, hence deformations of the sheaf $\sF\vert_{U_{\beta}}$ by taking the sheaf associated to the $0$-th cohomology. In fact, the condition $(d_{\sQ_{\beta}}+l_{\beta})^2=0$ is equivalent to require $d_{\mathfrak{L}_0}l_{\beta}+\frac{1}{2}[l_{\beta},l_{\beta}]=0$, while the flatness follows from~\cite[Theorem A.31]{Ser} since every cofibrant complex is degreewise projective, see e.g.~\cite[Lemma 2.3.6]{Hov99}. Notice that for every $\alpha\leq\beta\leq\gamma$ we have a quasi-isomorphism
\[ (\sQ_{\beta}\otimes A,d_{\sQ_{\beta}}+l_{\beta})\otimes_{(A_{\beta}\otimes A)}(A_{\gamma}\otimes A) \to (\sQ_{\gamma}\otimes A,d_{\sQ_{\gamma}}+l_{\gamma}) \]
so that the induced map between deformations
\[ \xymatrix{ H^0(\sQ_{\beta}\otimes A,d_{\sQ_{\beta}}+l_{\beta})\otimes_{(A_{\beta}\otimes A)}(A_{\gamma}\otimes A) \ar@{->}[rr]^-{\cong} \ar@{->}[dr] & & H^0(\sQ_{\gamma}\otimes A,d_{\sQ_{\gamma}}+l_{\gamma}) \ar@{->}[dl] \\
 & \sF(U_{\beta}) & } \]
is an isomorphism. This means that a Maurer-Cartan element $l^{\alpha}=\{l_{\beta}\}_{\beta\geq\alpha}\in\Hom_{A_{\pal}}^{1}(\sQ,\sQ)_{\alpha}\otimes\mA$ is essentially a deformation of the sheaf $\sF\vert_{U_{\alpha}}$.

Now consider a Maurer-Cartan element $l=\{l^{\alpha}\}_{\alpha\in\sN_0}\in\prod\limits_{\alpha\in\sN_0}\Hom_{A_{\pal}}^{1}(\sQ,\sQ)_{\alpha}\otimes\mA$, so that each $l^{\alpha}$ is a Maurer-Cartan element in $\Hom_{A_{\pal}}^{1}(\sQ,\sQ)_{\alpha}\otimes\mA$. In order to glue the deformations associated to each $l^{\alpha}$, we need to require the existence of an isomorphism
\[ (\sQ_{\beta}\otimes A,d_{\sQ_{\beta}}+l_{\beta}^{\alpha})\otimes_{(A_{\beta}\otimes A)}(A_{\gamma}\otimes A) \xrightarrow{f} (\sQ_{\beta}\otimes A,d_{\sQ_{\beta}}+l_{\beta}^{\alpha'})\otimes_{(A_{\beta}\otimes A)}(A_{\gamma}\otimes A) \]
lifting the identity for every $\alpha,\alpha'\in\sN_0$ and every $\beta\in\bar{\sN}$ such that $\alpha,\alpha'\leq\beta$. Since $f$ lifts the identity on $\sQ_{\beta}$, then $f=e^{m^{(\alpha,\alpha')}_{\beta}}$ for some $m_{\beta}^{(\alpha,\alpha')}\in\Hom_{A_{\pal}}^0(\sQ_{\beta},\sQ_{\beta})\otimes\mA$. The commutativity with the differential is equivalent to the relation $d_{\sQ_{\beta}}+l_{\beta}^{\alpha} =  e^{m_{\beta}^{(\alpha,\alpha')}} (d_{\sQ_{\beta}}+l_{\beta}^{\alpha'}) e^{-m_{\beta}^{(\alpha,\alpha')}}$, i.e. $l_{\beta}^{\alpha'} = e^{m_{\beta}^{(\alpha,\alpha')}}\ast l_{\beta}^{\alpha}$.
Therefore for every $(\alpha,\alpha')\in\bar{\sN}_1$ all these isomorphisms are collected by the element $(\alpha,\alpha')\in \Hom_{A_{\pal}}^0(\sQ,\sQ)_{\alpha\cup\alpha'}\otimes\mA$.

Observe that in order to satisfy the cocycle condition on the $0$-th cohomology, we need to require that for every $(\alpha,\alpha',\alpha'')\in\bar{\sN}_2$ there exists an element $n^{(\alpha,\alpha',\alpha'')}\in\Hom_{A_{\pal}}^{-1}(\sQ,\sQ)_{\alpha\cup\alpha'\cup\alpha''}$ such that
\[ m^{(\alpha',\alpha'')}_{\gamma} \bullet (-m^{(\alpha,\alpha'')}_{\gamma}) \bullet m^{(\alpha,\alpha')}_{\gamma} = \left[d+l^{\alpha'}_{\gamma},n_{\gamma}^{(\alpha,\alpha',\alpha'')}\right] \]
for every $\gamma\geq(\alpha,\alpha',\alpha'')$.

Summing up all the above discussion, we have a natural transformation defined for every $A\in\Art_{\K}$ by
\[ \varphi_A\colon H^1_{\mathfrak{L}}(A) \longrightarrow \Def_{\sF}(A) \; , \qquad \qquad \left(\{l^{\alpha}\}_{\alpha\in\sN_0},\{m^{(\alpha,\alpha')}\}_{(\alpha,\alpha')\in\bar{\sN}_1}\right) \mapsto (\sF_A\to\sF) \]
where $\sF_A$ is the sheaf obtained gluing together the deformations associated to each $l^{\alpha}$ through the isomorphisms $e^{m^{(\alpha,\alpha')}}$.

\begin{proposition}\label{prop.FIM}
The natural transformation $\varphi\colon\Def_{\sF}\to H^1_{\mathfrak{L}}$ defined above is a natural isomorphism.
\begin{proof}
For simplicity we assume the replacement $\sQ$ to belong to $\PM^{\leq 0}(A_{\pal})$, i.e. $\sQ_{\alpha}$ is concentrated in non-positive degrees for every $\alpha\in\sN$. Notice that by Remark~\ref{rmk.boundedpseudomodule} such a replacement always exists, and our assumption is not restrictive since for every pair of cofibrant replacements $\sQ\to\Upsilon^{\ast}\sF\leftarrow\sQ'$ the DG-Lie algebras $\End_{A_{\pal}}^{\ast}(\sQ)$ and $\End_{A_{\pal}}^{\ast}(\sQ')$ are quasi-isomorphic.

In order to prove the claim, fix $A\in\Art_{\K}$ and take an isomorphism between deformations $f\colon\sF_A$ and $\sF'_A$, associated to $\left(\{l^{\alpha}\}_{\alpha\in\sN_0},\{m^{(\alpha,\alpha')}\}_{(\alpha,\alpha')\in\bar{\sN}_1}\right)$ and $\left(\{\lambda^{\alpha}\}_{\alpha\in\sN_0},\{\mu^{(\alpha,\alpha')}\}_{(\alpha,\alpha')\in\bar{\sN}_1}\right)$ respectively. For every $\alpha\in\sN_0$ and every $\beta\geq\alpha$, the restriction of $f$ to each $U_{\alpha}$ lifts to isomorphisms
\[ (\sQ_{\beta}\otimes A,d_{\sQ_{\beta}}+l^{\alpha}_{\beta})\to(\sQ_{\beta}\otimes A,d_{\sQ_{\beta}}+\lambda^{\alpha}_{\beta}) \]
that reduce to the identity modulo the maximal ideal $\mA$. Therefore all these isomorphisms are of the form $e^{a^{\alpha}_{\beta}}$ for some $\{a^{\alpha}\}\in\prod\limits_{\alpha\in\sN_0}\Hom^0_{A_{\pal}}(\sQ,\sQ)_{\alpha}\otimes\mA$. Again, the commutativity with the differentials is equivalent to the relations
\[ e^{a^{\alpha}_{\beta}}\ast l_{\beta}^{\alpha} = \lambda_{\beta}^{\alpha} \; ,\qquad \qquad \text{ for every } \beta\geq\alpha \; . \]

We are only left with the proof that $\varphi_A$ is surjective for every $A\in\Art_{\K}$. To this aim, take a deformation $\sF_A\to\sF$ in $\Def_{\sF}$ and fix $\alpha\in\sN_0$. Notice that for every $\beta\geq\alpha$ in $\sN$ the map $\sQ_{\beta}\to\sF(U_{\beta})$ lifts to surjective quasi-isomorphisms $(\sQ_{\beta}\otimes A, d+l_{\beta}^{\alpha}) \to \sF_A(U_{\beta})$  of DG-modules over $A_{\beta}\otimes A$, for some $l^{\alpha}\in\Hom_{A_{\pal}}^1(\sQ,\sQ)_{\alpha}\otimes\mA$. The gluing data correspond to elements $m^{(\alpha,\alpha')}\in\Hom_{A_{\pal}}^0(\sQ,\sQ)_{\alpha\cup\alpha'}\otimes\mA$ for every $(\alpha,\alpha')\in\bar{\sN}_1$; moreover, for every $\beta\geq\alpha\cup\alpha'$ each isomorphism $e^{m^{(\alpha,\alpha')}_{\beta}}$ lifts the identity in the $0$-th cohomology, and liftings are unique up to homotopy.
\end{proof}
\end{proposition}

The argument used in Proposition~\ref{prop.FIM} is similar to the Kodaira-Spencer approach to deformations of a locally free sheaf $\sE$ of $\Oh_X$-modules on a complex manifold,~\cite{KS}, and in fact closely follows the one given in~\cite{FIM} to show that deformations of a quasi-coherent sheaf $\sF$ are controlled by the sheaf of DG-Lie algebras $\sE nd^{\ast}(\sE)$ for any given locally free resolution $\sE\to\sF$. The main advantage of our approach relies on the fact that we do not assume the existence of such a resolution.

\begin{theorem}\label{thm.FIM}
Let $X$ be a finite dimensional Noetherian separated scheme over $\K$, and let $\sF$ be a quasi-coherent sheaf on it. Fix a cofibrant replacement $\sQ\to\Upsilon^{\ast}\sF$. Then there exists a natural isomorphism $\Def_{\Tot_{TW}(\mathfrak{L})}\longrightarrow \Def_{\sF}$, where $\mathfrak{L}$ is the semicosimplicial DG-Lie algebra associated to $\sQ$, see Definition~\ref{def.semicosimplicialL}.

Hence by Remark~\ref{rmk.diagram} we have natural isomorphisms $\Def_{\End_{A_{\pal}}^{\ast}(\sQ)}\cong\Def_{\Tot_{TW}(\mathfrak{L})}\cong\Def_{\sF}$.
\begin{proof}
It has been already observed in Remark~\ref{rmk.diagram} that $\Def_{\End_{A_{\pal}}^{\ast}(\sQ)}\cong\Def_{\Tot_{TW}(\mathfrak{L})}$. Therefore, by Lemma~\ref{lemma.positivecohomology} and~\cite[Theorem 7.6]{FIM}, it is sufficient to prove that $\Def_{\sF}=H^1_{\mathfrak{L}}$. The statement now follows by Proposition~\ref{prop.FIM}.
\end{proof}
\end{theorem}

In particular, by Corollary~\ref{corollary.diagram} we recover the well-known fact that $T^1\Def_{\sF}=\Ext^1(\sF,\sF)$ and obstructions are contained in $\Ext^2(\sF,\sF)$.

\bigskip

\subsection{Deformations via $A_{\pal}$-modules}\label{section.deformationpseudomod}

In this subsection we present another proof of Theorem~\ref{thm.FIM} without using semicosimplicial techniques.

\begin{theorem}\label{thm.pseudoDef}
Let $X$ be a finite dimensional Noetherian separated scheme over $\K$, and let $\sF$ be a quasi-coherent sheaf on it. Fix a cofibrant replacement $\sQ\to\Upsilon^{\ast}\sF$. Then there exists a natural isomorphism $\Def_{\End^{\ast}_{A_{\pal}}(\sQ)} \longrightarrow \Def_{\sF}$.

Hence by Remark~\ref{rmk.diagram} we have natural isomorphisms $\Def_{\Tot_{TW}(\mathfrak{L})}\cong\Def_{\End_{A_{\pal}}^{\ast}(\sQ)}\cong\Def_{\sF}$.
\begin{proof}
For simplicity we assume the replacement $\sQ$ to belong to $\PM^{\leq 0}(A_{\pal})$, i.e. $\sQ_{\alpha}$ is concentrated in non-positive degrees for every $\alpha\in\sN$. Notice that by Remark~\ref{rmk.boundedpseudomodule} such a replacement always exists, and our assumption is not restrictive since for every pair of cofibrant replacements $\sQ\to\Upsilon^{\ast}\sF\leftarrow\sQ'$ the DG-Lie algebras $\End_{A_{\pal}}^{\ast}(\sQ)$ and $\End_{A_{\pal}}^{\ast}(\sQ')$ are quasi-isomorphic, hence inducing isomorphic deformation functors $\Def_{\End_{A_{\pal}}^{\ast}(\sQ)}\cong\Def_{\End_{A_{\pal}}^{\ast}(\sQ')}$.

Our first goal is to explicitly define a natural transformation $\varphi\colon\Def_{\End^{\ast}_{A_{\pal}}(\sQ)} \longrightarrow \Def_{\sF}$. To every object $\eta=\{\eta_{\alpha}\}_{\alpha\in\sN}\in\MC\left(\Hom^{\ast}_{A_{\pal}}(\sQ,\sQ)\otimes A\right)$ there are associated (local) deformations
\[ H^0(\sQ_{\alpha}\otimes A,d_{\sQ_{\alpha}}+\eta_{\alpha}) \to \sF(U_{\alpha})\; , \qquad \qquad \alpha\in\sN \]
where each $H^0(\sQ_{\alpha}\otimes A,d_{\sQ_{\alpha}}+\eta_{\alpha})$ is $A$-flat by~\cite[Theorem A.31]{Ser}. Here the Maurer-Cartan equation is equivalent to the condition $(d_{\sQ_{\alpha}}+\eta_{\alpha})^2=0$. Moreover, for every $\alpha\leq\beta$ the map
\[ H^0(\sQ_{\alpha}\otimes A,d_{\sQ_{\alpha}}+\eta_{\alpha})\otimes_{(A_{\alpha}\otimes A)}(A_{\beta}\otimes A) \to H^0(\sQ_{\beta}\otimes A,d_{\sQ_{\beta}}+\eta_{\beta}) \]
is an isomorphism because $\sQ$ is quasi-coherent in $\PM(A_{\pal})$ by Remark~\ref{rmk.qcoh-homotopyinvariance}. Now, for every $\alpha\leq\beta\leq\gamma$ there is a commutative diagram
\[ \xymatrix{	\sQ_{\alpha}\otimes_{A_{\alpha}}A_{\gamma} \ar@{->}[rr]_{q_{\alpha\beta}\otimes \Id_{A_{\gamma}}} \ar@{->}@/^2pc/[rrr]^{q_{\alpha\gamma}} & & \sQ_{\beta}\otimes_{A_{\beta}}A_{\gamma} \ar@{->}[r]_-{q_{\beta\gamma}} & \sQ_{\gamma} 	} \]
inducing the cocycle conditions on the deformations $\{H^0(\sQ_{\alpha}\otimes A,d_{\sQ_{\alpha}}+\eta_{\alpha}) \to \sF(U_{\alpha})\}_{\alpha\in\sN}$. Hence they glue together in a global deformation $\sF_A^{\eta}\to \sF$, with $\sF_A$ flat over $\Spec(A)$. Define the natural transformation $\varphi\colon\Def_{\End^{\ast}_{A_{\pal}}(\sQ)} \longrightarrow \Def_{\sF}$ as $\varphi_A\colon \eta \mapsto (\sF_A^{\eta}\to \sF)$ on every $A\in\Art_{\K}$. In order to show that $\varphi$ is well-defined, take two Maurer-Cartan elements $\eta,\xi\in\Hom^1_{A_{\pal}}(\sQ,\sQ)\otimes\mA$ and suppose that there exists an element $a=\{a_{\alpha}\}_{\alpha\in\sN}\in\Hom^0_{A_{\pal}}(\sQ,\sQ)\otimes\mA$ such that $e^a\ast\eta=\xi$. The last condition is equivalent to require that the maps in the square
\[ \xymatrix{	(\sQ_{\alpha}\otimes A, d_{\sQ_{\alpha}}+\eta_{\alpha}) \otimes_{(A_{\alpha}\otimes A)}(A_{\beta}\otimes A) \ar@{->}[rr]^{e^{a_{\alpha}}\otimes\Id_{(A_{\beta}\otimes A)}} \ar@{->}[d] & & (\sQ_{\alpha}\otimes A, d_{\sQ_{\alpha}}+\xi_{\alpha}) \otimes_{(A_{\alpha}\otimes A)}(A_{\beta}\otimes A) \ar@{->}[d] \\
(\sQ_{\beta}\otimes A, d_{\sQ_{\beta}}+\eta_{\beta}) \ar@{->}[rr]^{e_{\beta}} & & (\sQ_{\beta}\otimes A, d_{\sQ_{\beta}}+\xi_{\beta}) } \]
commute with differentials for every $\alpha\leq\beta$ in $\sN$. Therefore the associated deformations $\sF_A^{\eta}\to\sF$ and $\sF_A^{\xi}\to\sF$ are isomorphic.

We are left with the proof that $\varphi$ is a natural isomorphism. Fix $A\in\Art_{\K}$ and take an isomorphism between deformations $f\colon\sF_A^{\eta}$ and $\sF_A^{\xi}$, associated to $\eta=\{\eta_{\alpha}\}_{\alpha\in\sN}$ and $\xi=\{\xi_{\alpha}\}_{\alpha\in\sN}$ respectively. For every $\alpha\leq\beta$, the restriction of $f$ to each $U_{\alpha}$ lifts to isomorphisms
\[ (\sQ_{\alpha}\otimes A,d_{\sQ_{\alpha}}+\eta_{\alpha})\to (\sQ_{\alpha}\otimes A,d_{\sQ_{\alpha}}+\xi_{\alpha}) \]
that reduce to the identity modulo the maximal ideal $\mA$.
Therefore all these isomorphisms are of the form $e^{a_{\alpha}}$ for some $a=\{a_{\alpha}\}_{\alpha\in\sN}\in\Hom^0_{A_{\pal}}(\sQ,\sQ)\otimes\mA$. As above, the commutativity with the differentials is equivalent to the relations $e^{a}\ast \eta = \xi$, so that $\varphi_A$ is injective.

In order to show that $\varphi$ is surjective, fix $A\in\Art_{\K}$ and take a deformation $\sF_A\to\sF$ in $\Def_{\sF}$. Notice that for every $\alpha$ in $\sN$ the map $\sQ_{\alpha}\to\sF(U_{\alpha})$ lifts to surjective quasi-isomorphisms $(\sQ_{\alpha}\otimes A, d+\eta_{\alpha}) \to \sF_A(U_{\alpha})$ of DG-modules over $A_{\alpha}\otimes A$, for some $\eta_{\alpha}\in\Hom_{A_{\pal}}^1(\sQ,\sQ)\otimes\mA$.
\end{proof}
\end{theorem}

In particular, by Theorem~\ref{thm.REndpseudomodules} we recover the well-known fact that $T^1\Def_{\sF}=\Ext^1(\sF,\sF)$ and obstructions are contained in $\Ext^2(\sF,\sF)$.

If the sheaf $\sF$ admits a locally free resolution $\sE\to\sF$ then there exists a natural isomorphism of deformation functors $\Def_{\Tot_{TW}(\mathfrak{h})}\cong\Def_{\sF}$ by Theorem~\ref{thm.locallyfreeDef} and Theorem~\ref{thm.pseudoDef}.

\end{document}